\documentclass[preprint, 3p]{elsarticle}

\usepackage{lipsum}
\makeatletter
\def\ps@pprintTitle{%
\let\@oddhead\@empty
\let\@evenhead\@empty
\def\@oddfoot{}%
\let\@evenfoot\@oddfoot}
\makeatother

\usepackage{algorithm}
\usepackage{algpseudocode}
\usepackage{amsmath}
\usepackage{amssymb}
\usepackage{amsthm}
\usepackage{amscd}
\usepackage{appendix}
\usepackage{bbm}
\usepackage{braket}
\usepackage{bookmark}
\usepackage{booktabs}
\usepackage{caption}
\usepackage{cases}
\usepackage{enumitem}
\usepackage[T1]{fontenc} %codifica del font di uscita
\usepackage{float}
\usepackage{lineno,hyperref}
\modulolinenumbers[5]
\usepackage{indentfirst}
\usepackage[utf8]{inputenc}
\usepackage[autostyle,italian=guillemets]{csquotes}
\usepackage{graphicx}
\graphicspath{{Immagini/}}
\DeclareGraphicsExtensions{{.png},{.pdf},{.jpg}}
\graphicspath{{./}}
\usepackage{lmodern}
\usepackage{lipsum} %genera testo fittizio
\usepackage{listings}
\usepackage{yhmath}
\usepackage{mathtools}
\usepackage{mathrsfs}

\usepackage{microtype} %Sistema il carattere in modo che sia pi bello da vedere nella riga
\usepackage{multirow}
\usepackage{setspace}
\usepackage{standalone}
\usepackage{siunitx}
\usepackage{subcaption}
\usepackage{tabularx}
\usepackage[table]{xcolor}
\usepackage{wasysym}
\usepackage{wrapfig}
\usepackage{url} %permette di scrivere indirizzi internet

\usepackage[overload]{empheq}
\usepackage{cleveref}

\usepackage{pgfplots}
\pgfplotsset{compat=newest}
\usetikzlibrary{plotmarks}
\usetikzlibrary{arrows.meta}
\usepgfplotslibrary{patchplots}
\usepackage{grffile}

\newdefinition{rmk}{Remark}
\newproof{prf}{Proof}

\begin{document}
	
\begin{frontmatter}
\title{A staggered-in-time and non-conforming-in-space numerical framework for realistic cardiac electrophysiology outputs}

\author[1]{Elena Zappon\corref{cor1}}
\ead{elena.zappon@polimi.com}
\author[1]{Andrea Manzoni}
\ead{andrea1.manzoni@polimi.com}
\author[1,2]{Alfio Quarteroni}
\ead{alfio.quarteroni@polimi.it}

\cortext[cor1]{Corresponding author}
		
\address[1]{MOX - Dipartimento di Matematica, Politecnico di Milano, P.zza Leonardo da Vinci 32, I-20133 Milano, Italy}
	
\address[2]{Institute of Mathematics, Ecole Polytechnique Federale de Lausanne, Station 8,  CH-1015 Lausanne, Switzerland (Professor Emeritus)}

\begin{abstract}
Computer-based simulations of non-invasive cardiac electrical outputs, such as electrocardiograms and body surface potential maps, usually entail severe computational costs due to the need of capturing fine-scale processes and to the complexity of the heart-torso morphology. In this work,
we model cardiac electrical outputs by employing a coupled model consisting of a reaction-diffusion model - either the bidomain model or the most efficient pseudo-bidomain model - on the heart, and an elliptic model in the torso. We then solve the coupled problem with a segregated and staggered in-time numerical scheme, that allows for independent and infrequent solution in the torso region. To further reduce the computational load, main novelty of this work is in introduction of an interpolation method at the interface between the heart and torso domains, enabling the use of non-conforming meshes, and the numerical framework application to realistic cardiac and torso geometries. The reliability and efficiency of the proposed scheme is tested against the corresponding state-of-the-art bidomain-torso model. Furthermore, we explore the impact of torso spatial discretization and geometrical non-conformity on the model solution and the corresponding clinical outputs. The investigation of the interface interpolation method provides insights into the influence of torso spatial discretization and of the geometrical non-conformity on the simulation results and their clinical relevance.
\end{abstract}
\end{frontmatter}

\section{Introduction}
In the last decades computer-based simulations of non-invasive cardiac electrical outputs, such as  electrocardiograms (ECGs) and body surface potential maps (BSPMs), have been widely used to enhance the understanding of cardiac diseases - \emph{e.g.} cardiac ischemia \cite{GEMMELL2020103895,loewe2018cardiac,loewe2015ecg,wilhelms2011comparing}, ventricular tachycardia \cite{lopez2019personalized,kania2017prediction,green1995clinical}, and atrial fibrillation \cite{pezzuto2018beat,rodrigo2017highest,zhou2016noninvasive} . Thanks to the use of high fidelity, patient-specific cardiac emulators \cite{GILLETTE2021102080,Gillette2018,MONACI2020104005,VILLONGCO2014305}, therapy planning \cite{green1995clinical,Monaci2022,Monaci2021Auto}, and drug actions \cite{noble2008computational} have also been remarkably improved. The main factors behind the success of cardiac mathematical simulations are the accurate representation of the cardiac and torso anatomies, and in an exhaustive physically-detailed model able to describe the undergoing biological processes. 

The effort to bridge the gap between real-world (in vivo or in vitro) observations and computational (in silico) simulations \cite{lamata2014images,vadakkumpadan2012image,2010Bishop,Zheng2007,yushkevich2006user} has led to the development of detailed and structurally-complex heart-torso geometrical representations. The preservation of accurate anatomical models, coupled with complex mathematical representations of cardiac physiology, requires however a high spatial resolution for the discretization of the heart-torso geometry \cite{arevalo2016arrhythmia,Niederer2011bench}, and this level of detail comes with a significant increase of computational costs.

A mathematical model for the generation of electrophysiology (EP) outputs is usually obtained by coupling a model for the simulation of the electrical current in the heart with a model of the human body - here referred as \emph{torso} - treated as a passive conductor \cite{Pullan2005,Lines2003,Malmivuo1995}. 
Realistic simulations of the complex cardiac behavior require a high number of interrelated features of the heart model \cite{FEDELE2023115983,Gerach2021,Gillette2021purkinje}, and determine an exponential growth of the computational costs. The development of efficient diffusive models in the torso  \cite{ColliFranzone2014,Boulakia2010,Potse2009,Bradley2000,Aoki1987}, such as the lead field method \cite{multerer2021uncertainty,Potse2018,mcfee1953electrocardiographic}, or the boundary element method (BEM) \cite{Potse2009,Schuler2019,fischer2000bidomain}, has become paramount to contained 
such computational effort. However, both torso representations show severe limitations in reproducing certain cardiac electrophysiology phenomena or specific clinical outputs. Indeed, the efficiency of the lead field methods relies on the simplified monodomain model \cite{Pullan2005,ColliFranzone2014,Quarteroni2019,Sundnes2006} for the EP simulation, and is therefore not suited to represents EP events influenced by the potential field in the extracellular space. Moreover, the lead field outputs are one-dimension EP signals, \emph{e.g.} the 12-lead ECG, whereas body surface potential maps cannot be computed. On the other hand, boundary BEM-based approaches lack the capability to incorporate cardiac muscle anisotropy in the extracellular space, limiting their applicability to pathological conditions \cite{Schuler2019}. Additionally, these methods cannot efficiently handle more complex descriptions of cardiac function, such as the coupling of electromechanical or electro-mechano-fluid models, as they rely on an static representation of the torso geometry.

The most biophysically detailed approach for generating cardiac EP outputs is represented by the \emph{fully-coupled heart-torso model} (FCHT), obtained by a perfect coupling - \emph{i.e.} prescribing interface continuity of the extracellular field at the heart-torso interface - of the bidomain model \cite{Pullan2005,ColliFranzone2014,Quarteroni2019} for the heart activation, and a Laplace diffusion equation for the human torso \cite{Pullan2005,Lines2003,Potse2006}. The FCHT solution is the electric field evolution in time on both heart and torso anatomies, and ECGs and BSPMs can be conveniently obtained by a suitable post-processing  of the simulation outputs. The efficient solution of the FCHT has been investigated by using either monolithic \cite{Lines2003,Pullan2005,Sundnes2002} or splitting schemes \cite{Lines2003,Boulakia2010,SUNDNES200155}, as well as matrix transfer operators \cite{Boulakia2010}, but remains a challenging task. Indeed, solving the bidomain model alone is 10 times more expensive than solving the monodomain model \cite{plank2009generation}, whereas the overall increase of number of degrees of freedom (DoFs) due to the torso geometry, as well as the additional costs of solving the diffusive model in the torso, make the computational expenses often unfeasible \cite{GILLETTE2021102080,Boulakia2010,Bishop2011_2}, unless coarse discretizations are used. 

An alternative that combines the monodomain model with the detailed solution of the FCHT model is reported in \cite{Potse2006}, where the cardiac EP dynamic is simulated with the monodomain model employing a fine time discretization, whereas the extracellular field is computed with a coarse time discretization (i) on the cardiac domain by solving the elliptic portion of the bidomain model, and (ii) on the torso domain by solving an elliptic problem. In \cite{Bishop2011_2,Bishop2011}, the authors demonstrate that, by applying an augmentation approach and carefully choosing the conduction velocity, the monodomain model can replicate the solution profiles of the bidomain one, preserving the computational efficiency. The reliability of the coupled monodomain-elliptic model, referred to as the \emph{pseudo-bidomain} model, has been previously demonstrated in \cite{Boulakia2010}. The authors showed that while the choice of the EP model influences the amplitude of the electric signals, both the pseudo-bidomain and FCHT approaches accurately capture the overall shape of the ECG. This finding highlights the comparable performance of the pseudo-bidomain model in reproducing key features of the ECG, providing evidence of its reliability as an alternative to the FCHT model.

In this study, we employ a simplified version of the FCHT model. The cardiac domain is coupled with either the bidomain or the pseudo-bidomain model, while the torso domain is represented using an elliptic model. At the interface between the heart and torso, we enforce the continuity of the extracellular solution by applying Dirichlet interface conditions on the torso model. However, we neglect the continuity of fluxes and assume complete isolation of the heart from the surrounding torso. This is achieved by imposing homogeneous Neumann interface conditions on the cardiac electrophysiology model. We then implement a segregated and staggered in-time numerical scheme, which enable for an independent solution of the electrophysiology problem on the cardiac domain, and a subsequent and infrequent solve of the elliptic model in the torso one. The main novelties of this work, however, lie in (i) the introduction of an interpolation method at the heart-torso interface, (ii) the application of the proposed staggered in-time and non-conforming in-space numerical framework on realistic and detailed cardiac and torso geometries, and (iii) on the investigation of the effect of the torso geometry and its spatial discretization on clinical electrophysiology output. Indeed, the interpolation method at the domains interface enables the independent treatment of the spatial discretization of the cardiac and torso domains, allowing for \emph{mesh-based non-conformity}. Furthermore, it allows the use of a heart domain that may have a different shape or position compared to the portion of cardiac tissue within the torso domain, thus enabling \emph{geometrical non-conformity}. By incorporating this novel approach, the flexibility of the model is significantly enhanced, allowing for variations in the spatial discretization and geometry of both the heart and torso domains. Consequently, the model captures the spatial and temporal evolution of the extracellular field within these domains. Furthermore, important clinical outputs such as ECGs and BSPMs can be efficiently extracted from the model results during a post-processing stage. This comprehensive framework enables the simulation and analysis of complex cardiac electrophysiology phenomena, providing valuable insights into the behavior of the extracellular field and its impact on clinical measurements.

We evaluate the accuracy of our approach in reproducing 12-lead ECG traces and BSPMs compared to the fully-coupled heart-torso (FCHT) model. We also highlight the computational benefits of our method. By leveraging the interpolation method at the interface between the heart and torso domains, we investigate the influence of torso discretization on the torso solution and the resulting ECG signals, demonstrating that the mesh size of the torso has a limited impact on the overall problem solution. Furthermore, we examine the performance of our framework on heart-torso domains with geometrical non-conformity. By imposing rigid transformations on the heart and projecting the EP solution onto a common reference torso domain, we assess the quantitative reliability of our model under small non-conformities and the qualitative fidelity under large non-conformities. Overall, our findings demonstrate the effectiveness and versatility of our proposed approach in accurately simulating and analyzing cardiac electrophysiology phenomena in various scenarios.
 
The paper is organized as follows: in Section \ref{Sec:models} we provide a brief overview of the mathematical models involved and the numerical framework utilized; in Section \ref{Sec:numerical_results} we present the numerical results obtained using the proposed scheme; whereas in Section \ref{Sec:conclusions} we concludes the paper by summarizing the key findings and implications of the study.  

\section{Models and methods}
\label{Sec:models}
\begin{figure}[!t]
	\centering
	\includegraphics[width=0.65\textheight]{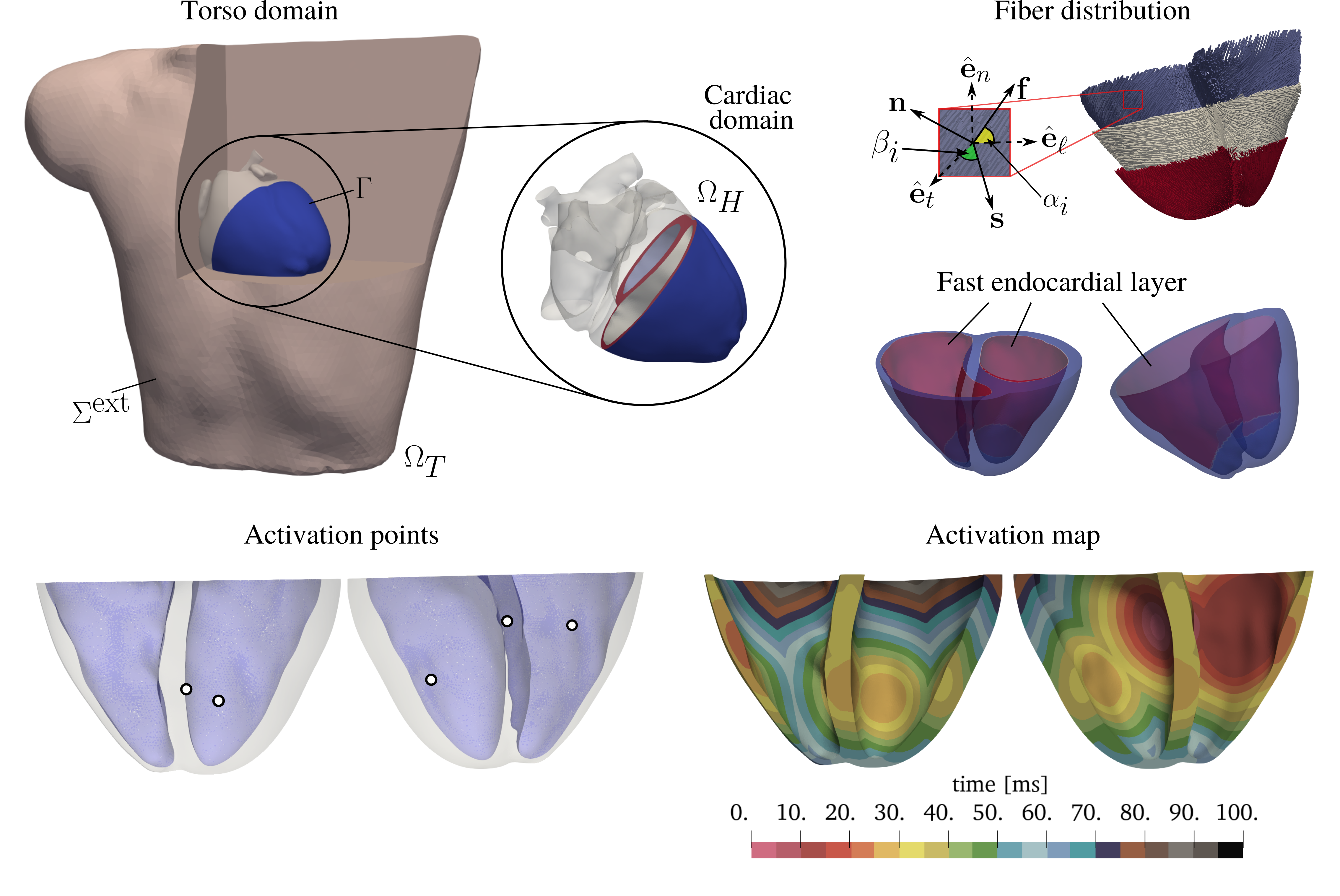}
	\caption{Top-left: section of the torso domain $\Omega_T$, and of the cardiac domain $\Omega_H$. Here $\Omega_H$ is represented by the biventricular geometry, obtained by cutting the heart with a planar surface well below the mitral valve. Top-right: fiber distribution and fast endocardial layer computed on $\Omega_H$. Bottom-left: location of the five activation points on the endocardial surface. Bottom-right: activation map on the biventricular geometry \cite{DURRER1970}.}
	\label{Fig:EP_data}
\end{figure}

Let $\Omega \subset \mathbb{R}^3$ be an open, bounded domain representing the human body, partitioned into two non-overlapping subdomains $\Omega_H$ and $\Omega_T$ defining the volume occupied by the heart, and the remaining part of the body, that we generally define as torso domain, respectively (see Fig. \ref{Fig:EP_data}). We denote the common interface as $\overline{\Omega_H} \cap \overline{\Omega_T} = \Gamma$. 

Let $t \in (0,T]$ denote the independent variable representing time, and $\mathbf{x}$ represent the spatial coordinate. In this study, we employ a coupled model, referred to as the heart-torso model, which incorporates two distinct physics: electrophysiology, to simulate the activation of the cardiac tissue	 \cite{Pullan2005,Lines2003,ColliFranzone2014,Quarteroni2019}, and a passive conduction, to represent the propagation of the cardiac electric signal through the rest of the body \cite{Boulakia2010,Boulakia2011}. The unknowns of the heart-torso model are:
\[\begin{array}{ll}
V_m: \Omega_H \times (0,T] \longrightarrow \mathbb{R} & \text{transmembrane potential},\\
u_e: \Omega_H \times (0,T] \longrightarrow \mathbb{R} & \text{extracellular potential},\\
\mathbf{w}: \Omega_H \times (0,T] \longrightarrow \mathbb{R}^{N^{w}_{ion}} &\text{ionic gating variables},\\
\mathbf{c}: \Omega_H \times (0,T] \longrightarrow \mathbb{R}^{N^{c}_{ion}} &\text{ionic concentrations},\\
u_T: \Omega_T \times (0,T] \longrightarrow \mathbb{R} & \text{extracellular potential in the torso}.
\end{array}\]

The following section provides a detailed description of the model used in this study, including the coupling approach and the numerical strategies employed.

\subsection{Bidomain model}
The subcellular ionic dynamics are the main responsible of the generation of an action potential that yields variations of the intracellular $u_i$ and extracellular $u_e$ electric fields, and ultimately of the transmembrane potential $V_m = u_i - u_e$. 

The biophysically most detailed description of the evolution of the transmembrane and the extracellular potential is the the bidomain model, a nonlinear reaction-diffusion system of partial differential equations (PDEs) \cite{ColliFranzone2014,TenTusscher2006}, that reads as follows:

\begin{subequations}
	\label{Eq:bidomain_model}
	\begin{empheq}[left={\empheqlbrace\,}]{align}
	&\chi\left [ C_m \frac{\partial V_m}{\partial t} + I^{ion}(V_m,\mathbf{w},\mathbf{c})\right] - \nabla \cdot(\mathbf{D}_i\nabla (V_m + u_e)) = I^{app}(t) \qquad \text{in }\Omega_H \times (0,T]  \label{Eq:bid_1}\\
	&-\nabla \cdot(\mathbf{D}_i\nabla V_m) - \nabla ((\mathbf{D}_i + \mathbf{D}_e)\nabla u_e) = 0 \qquad \qquad\qquad\qquad\quad\quad \:\:\:\:\, \text{in }\Omega_H \times (0,T] \label{Eq:bid_2}\\
	& (\mathbf{D}_i\nabla (V_m + u_e)) \cdot \mathbf{n}_H = 0 \qquad\qquad\qquad\qquad\qquad\quad \qquad\qquad\qquad\:\:\:\:\,\text{on }\partial\Omega_H \times (0,T] \label{Eq:bid_5}\\
	&V_m= V_{m,0}, \qquad\qquad \qquad \qquad\qquad \qquad\qquad\qquad \qquad \qquad \qquad \quad \:\:\:\:\,\text{at }t = 0,
	\end{empheq}
\end{subequations}
where $C_m$ is the capacitance per unit area, $\chi$ is the surface-to-volume ratio on the membrane, $T$ is the maximum time of simulation, and $\mathbf{n}_H$ is the outward unit normal vector. 

The term $I^{\mathrm{ion}}(V_m,\mathbf{w},\mathbf{c})$ defines the ionic current flowing through the cardiac tissue, and depends on the vector of gating variable $\mathbf{w}$, describing the fraction of open ionic channels for a single cell, and the ionic concentration vector $\mathbf{c}$ \cite{ColliFranzone2014}. The dynamics of $\mathbf{w}$ and $\mathbf{c}$ is described by the ten Tusscher and Panfilov (TTP06) \cite{TenTusscher2006} model,  expressed by the following system of ordinary differential equations (ODEs)
\begin{subequations}
\label{Eq:ionic_model}
\begin{empheq}[left={\empheqlbrace\,}]{align}
&\frac{d \mathbf{w}}{d t} = \mathbf{H}(V_m,\mathbf{w}) \qquad\qquad\qquad\qquad\qquad\qquad\qquad\qquad\qquad \qquad\quad \:\:\:\,\text{in }\Omega_H \times (0,T] \label{Eq:bid_3}\\
&\frac{d \mathbf{c}}{d t} = \mathbf{G}(V_m,\mathbf{w},\mathbf{c})\qquad\qquad\qquad\qquad\qquad\qquad \qquad\qquad\qquad\quad\quad\:\:\:\:\:\,\text{in }\Omega_H \times (0,T] \label{Eq:bid_4}\\
&\mathbf{w} = \mathbf{w}_0,~ \mathbf{c} = \mathbf{c}_0\qquad \qquad\qquad \qquad\qquad\qquad \qquad \qquad \qquad \qquad \quad \:\:\:\:\,\text{at }t = 0.
\end{empheq}
\end{subequations}

We encode the myocardial anisotropy in the tensors $\mathbf{D}_{i,e}$, describing the intracellular and extracellular conductive capacities as
\begin{equation} 
\label{Eq:conduction_tensors}
\mathbf{D}_{i,e} = \sigma_{\ell}^{i,e}(\mathbf{f}\otimes\mathbf{f}) + \sigma_t^{i,e}(\mathbf{s}\otimes\mathbf{s}) + \sigma_n^{i,e}(\mathbf{n}\otimes\mathbf{n}).
\end{equation}
In the above equation, the vectors $[\mathbf{f},\mathbf{s},\mathbf{n}]$ represent the fiber $\mathbf{f}$, sheet $\mathbf{s}$, and sheet-normal $\mathbf{n}$ directions, and define an orthotropic reference system for the fiber generation. For this task, we employ the Laplace-Dirichlet-Rule-Based-Methods described in \cite{Piersanti2021} (see Fig. \ref{Fig:EP_data}), with the following input angles:
\begin{align}
\label{Eq:angle_fibers}
\alpha_{epi,LV} = \alpha_{epi,RV} = -\ang{60},\quad \alpha_{endo,LV} = \alpha_{endo,RV} = +\ang{60},\\
\beta_{epi,LV} = \beta_{epi,RV} = +\ang{25},\quad \beta_{endo,LV}= \beta_{endo,RV}= -\ang{65}.\end{align}
Different conduction coefficients $\sigma_{\ell,t,n}^{i,e}$ are prescribed to ensure the physiological conduction velocity of $0.6m/s$, $0.4m/s$ and $0.2m/s$ in longitudinal, transversal and
normal direction to the fibers, according to \cite{Pullan2005,ColliFranzone2014}. 

Finally, the cardiac conduction system is surrogated by means of a pacing protocol of five activation points selected on the endocardium, and a thin fast endocardial layer (see Fig. \ref{Fig:EP_data}), according to the physiological findings \cite{DURRER1970,Draper1959,Myerburg1978}. Such activation system is then encoded in the time-dependent forcing term $I^{\mathrm{app}}(t)$, that effectively provide the initial electric stimulus (see Fig. \ref{Fig:EP_data} for the activation map).

\begin{rmk}
Note that the bidomain model \eqref{Eq:bidomain_model}, being a pure Neumann problem, features multiple solutions. To overcome this issue, between the available strategies \cite{ColliFranzone2014,Pullan2005}, we impose a zero mean extracellular field on the cardiac domain by requiring that 
\[ \int_{\Omega_H} u_e(t) d\mathbf{x} = 0.
\]
\end{rmk}

\subsection{Pseudo-bidomain model}
In many cases, the evolution of the transmembrane potential can be effectively described using a simplified and more efficient version of the bidomain model, known as the monodomain model \cite{Potse2006}, whose solution represent the evolution in time of $V_m$  on the cardiac domain. However, it is important to note that the signal propagating through the human body is the extracellular potential $u_e$ generated by the heart. Therefore, while the transmembrane potential $V_m$ can be computed using the monodomain model, the extracellular potential $u_e$ needs to be recovered by solving the elliptic portion of the bidomain model \eqref{Eq:bid_2} on the cardiac domain. We refer to this combined solution of the monodomain and elliptic models as the pseudo-bidomain model.

The monodomain model is then derived by assuming similar anisotropy of the intracellular and extracellular domain, \emph{i.e.} $\sigma_{\ell,t,n}^{e} = \lambda \sigma_{\ell,t,n}^{i}$, for a given scalar $\lambda$ \cite{ColliFranzone2014}. A single PDE can then be employed to describe the temporal and spatial evolution of the transmembrane potential $V_m$. This PDE is given by

\begin{equation}
	\label{Eq:monodomain_model}
	\begin{cases}
	\chi\left [ C_m \frac{\partial V_m}{\partial t} + I^{ion}(V_m,\mathbf{w},\mathbf{c})\right] - \nabla \cdot(\mathbf{D}_m\nabla V_m ) = I^{app}(t) &\text{in }\Omega_H \times (0,T]  \\
	(\mathbf{D}_m\nabla V_m) \cdot \mathbf{n}_H = 0 &\text{on }\partial\Omega_H \times (0,T]\\
	V_m= V_{m,0}, ~\mathbf{w} = \mathbf{w}_0,~ \mathbf{c} = \mathbf{c}_0 &\text{at }t = 0,
	\end{cases}
\end{equation}
where $\mathbf{D}_m$ is used to describe the tissue conductive properties. In \cite{ColliFranzone2014}, the authors show that $\mathbf{D}_m$ is directly related to the extracellular and intracellular tensors \eqref{Eq:conduction_tensors} as 
$$ \mathbf{D}_m = \frac{\mathbf{D}_e\mathbf{D}_i}{\mathbf{D}_e+\mathbf{D}_i}.$$

As for the bidomain model, we employ the same pacing protocol encoded in the function $I^{app}(t)$ and the fast endocardial layer (see Fig. \ref{Fig:EP_data}) to simulate the cardiac activation, whereas the system of cardiac fibers are built by means of the Laplace-Rule-Based-Methods of \cite{Piersanti2021}. Ions dynamics is then represented also in this case by the TTP06 model. The extracellular potential $u_e$ in the heart is finally computed by solving the following model
\begin{equation}
\label{Eq:pseudo_bid}
-\nabla\cdot((\mathbf{D}_i+\mathbf{D}_e)\nabla u_e) = \nabla \cdot (\mathbf{D}_i\nabla V_m) \quad\text{in }\Omega_H \times (0,T].
\end{equation}

\subsection{The human body as a passive conductor}
Regarding the torso, a null intracellular current outside of the cardiac tissue \cite{Pullan2005,Sundnes2006,Krassowska1994} is generally assumed, whereas the cardiac extracellular potential propagates through the remaining part of the human body, and dissipate on its surface. Indeed, the torso is regarded as a passive conductor, isolated from the surrounding space.  

A mathematical model to describe the torso as a passive conductive volume is a generalized Laplace problem \cite{Pullan2005,Boulakia2010}, that we denote as \emph{torso model}, of the form
\begin{subequations}
	\label{Eq:torso}
	\begin{empheq}[left={\empheqlbrace\,}]{align}
&-\nabla \cdot(\mathbf{D}_T \nabla u_T) = 0 \qquad \text{in }\Omega_T\\
&\mathbf{D}_T\nabla u_T \cdot \mathbf{n}_T = 0 \qquad \quad\:\:\, \text{on }\Sigma^{\text{ext}}, \label{Eq:torso_2}
\end{empheq}
\end{subequations}
where $u_T$ is the extra-cellular potential field in the torso, and $\mathbf{n}_T$ is the outward unit normal with respect to $\Omega_T$. The insulation assumption is expressed by means of equation \eqref{Eq:torso_2}, where $\Sigma^{\text{ext}}$ denotes the external surface of the torso domain (see Fig. \ref{Fig:EP_data}). 

Different organs in the torso, such as lungs, bones, or the blood pool, have different conduction properties, but are generally considered isotropic. Although studies have shown that conduction heterogeneity can have a minimal impact on the ECGs and BSPMs  \cite{Bradley2000,Keller2010}, these findings are not specific to any particular mathematical model. For the sake of simplicity, in this study, we simulate a homogeneous torso with a conduction coefficient represented by the tensor $\mathbf{D}_T$, taking into account the average effect of the aforementioned components.

\subsection{Heart-torso coupling conditions}
\label{Sub:coupling_model}
Usual coupling conditions  between the electrophysiology model of the heart and the elliptic model in the torso enforce perfect transmission of the extracellular potential at the heart-torso interface $\Gamma$, ensuring continuity of both flux and solution  \cite{Pullan2005,ColliFranzone2014,Boulakia2010}. However, in this study, we employ a weak coupling strategy, which offers good reliabilty and improve the efficiency of the problem solution (refer to Section \ref{Sec:numerical_results}).

To this end, we postulate the heart to work in a fully-isolated way from the remaining human body, which amounts to imposing a zero extracellular flux at the heart-torso interface, that is
\[
\boldsymbol{\sigma}_e \nabla u_e \cdot \mathbf{n}_H = 0 \quad \text{on }\Gamma \times (0,T].
\]
As a results, we can add the following homogeneous Neumann boundary conditions 
\begin{equation}
\label{Eq:Neumann_one_way}
(\mathbf{D}_i+\mathbf{D}_e)\nabla u_e\cdot \mathbf{n}_H = -\mathbf{D}_i\nabla V_m \cdot \mathbf{n}_H \quad\text{on }\Gamma \times (0,T],
\end{equation}
to either of the EP models \eqref{Eq:bidomain_model} or \eqref{Eq:pseudo_bid}, and the EP quantities $V_m$ and $u_e$ can be computed independently of $u_T$.
To get the heart-torso coupling, we however require the continuity of the extracellular solution on $\Gamma$, obtained by imposing the following Dirichlet interface condition to the torso model \eqref{Eq:torso}
\begin{equation}
\label{Eq:dirichlet_torso}
u_T = u_e \quad \text{on }\Gamma \times (0,T].
\end{equation}
Summarizing, the cardiac potential $u_e$ is transferred from $\Omega_H$ to $\Omega_T$, whereas no electric feedback is considered from $\Omega_T$ to $\Omega_H$. We therefore define two one-way coupled models depending on the EP model employed, \emph{i.e.} the bidomain-torso model \eqref{Eq:bidomain_model}-\eqref{Eq:torso} and the pseudo-bidomain-torso model \eqref{Eq:monodomain_model}-\eqref{Eq:pseudo_bid}-\eqref{Eq:torso}. 

\begin{rmk}
	Note that the interface conditions \eqref{Eq:Neumann_one_way}-\eqref{Eq:dirichlet_torso} describe a unidirectional feedback from the heart to the torso, while neglecting the impact of the torso solution on the cardiac extracellular potential. The effect of such assumption has been investigated in \cite{Bishop2011}, where the authors study the outcome of the neglected feedback only on the electric signal propagation on the cardiac domain, whereas the corresponding effect on the ECGs has been previously observed only on idealized geometries \cite{Boulakia2010}. In Section \ref{Sec:numerical_results}, we will rigorously address this problem, investigating the ECG and BSPMs variations by solving either the FCHT model, the bidomain-torso model and the pseudo-bidomain torso model, on realistic heart and torso geometries. 
\end{rmk}

\subsection{A staggered in-time and non-conforming in-space numerical framework}
\label{Subsec:framework}
\begin{figure}[!t]
	\centering
	\includegraphics[width=0.25\textheight]{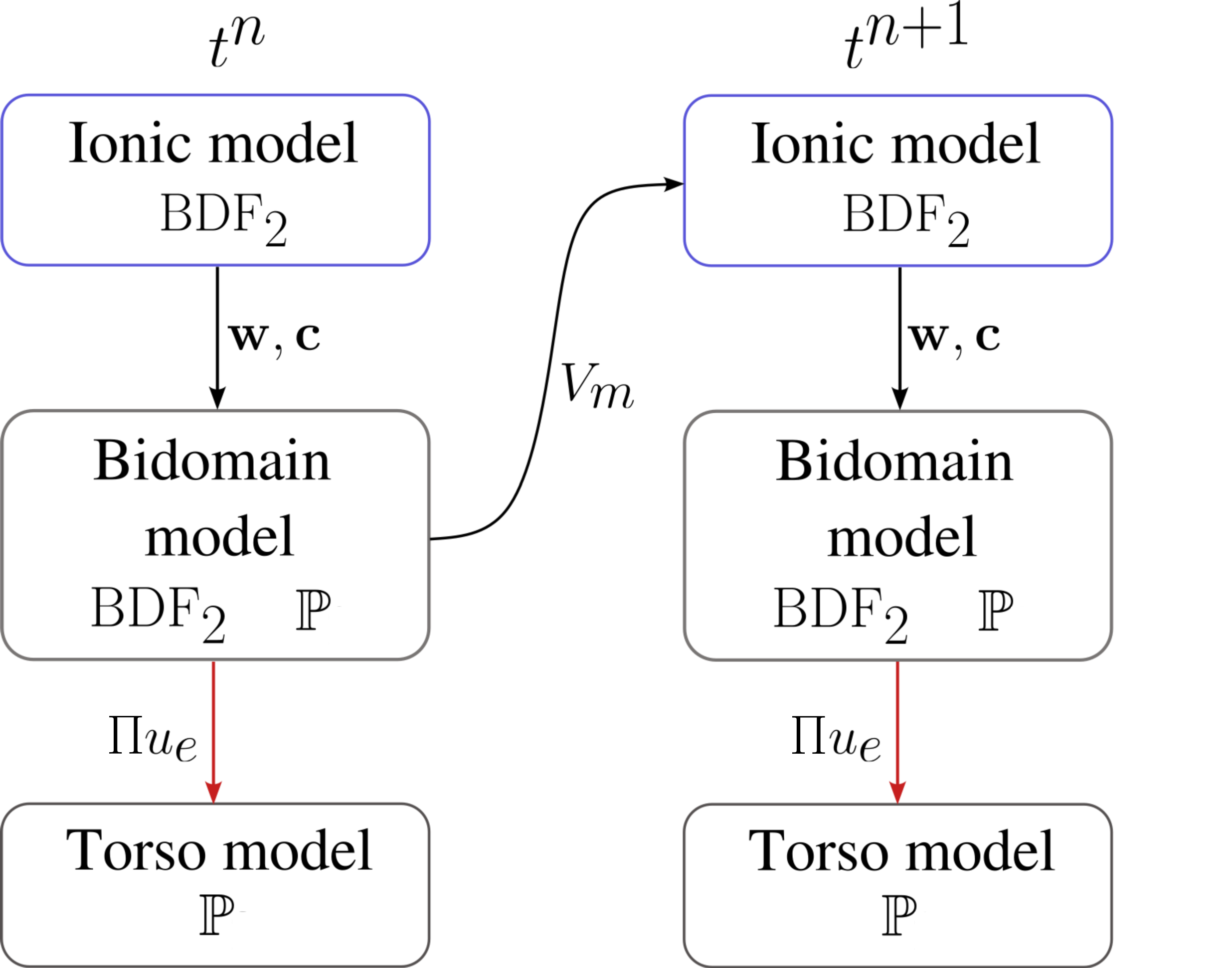}
	\quad
	\includegraphics[width=0.25\textheight]{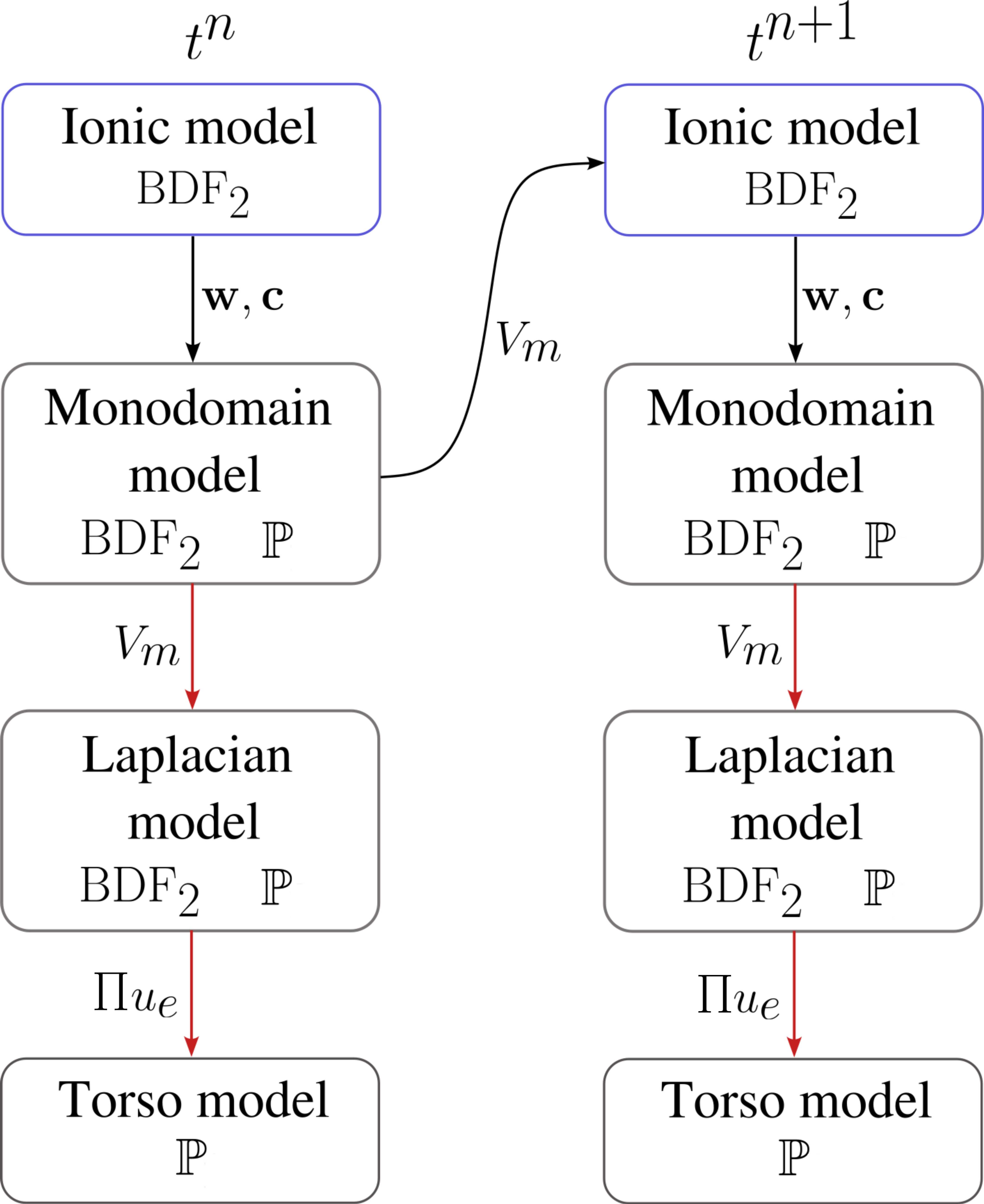}
	\caption{Schematic representation of the one-way coupled bidomain-torso model (left), and pseudo-bidomain-torso model (right).}
	\label{Fig:scheme_num}	
\end{figure}
The heart-torso coupled problem involves physical models with very distinct numerical requirements. The reaction-diffusion equations \eqref{Eq:bidomain_model}-\eqref{Eq:monodomain_model} governing the heart's electrical activity have fast dynamics and necessitate a fine temporal and spatial discretization on the cardiac domain $\Omega_H$ \cite{arevalo2016arrhythmia}. In contrast, the torso model \eqref{Eq:torso} has a smoother solution, and the clinical outputs of interest can be sampled at a lower frequency, allowing for a coarser spatial and temporal discretization on the torso domain $\Omega_T$.

To meet these requirements, we propose a segregated and staggered in-time computational framework that incorporates different time discretizations for the EP and torso models, using the weak coupling strategy described in Subsection \ref{Sub:coupling_model}.

Temporal discretization is then obtained by partitioning the time interval $(0,T]$ into sub-intervals $(t_i,t_{i+1}]$, $i=0,1,\dots,N_H$, $t_0=0$, and $t_{N_H} = T$, such that the step $\Delta t = t_{i+1} - t_i$, for all $i$, is small enough to fulfill the electrophysiology needs. Then, a multiple $m \Delta t$ of the EP time step is employed to solve the torso model. 

The time derivatives in the bidomain model \eqref{Eq:bidomain_model} and pseudo-bidomain model \eqref{Eq:pseudo_bid} are approximated using a second-order backward differentiation formula (BDF) scheme \cite{Ethier2008}. Additionally, an implicit-explicit (IMEX) scheme is employed for the ionic model \eqref{Eq:ionic_model}, with an implicit discretization for the gating variable and an explicit discretization for the ionic concentration \cite{Piersanti2022}. A semi-implicit discretization is used for the ionic term $I^{ion}$ in the reaction-diffusion equations \eqref{Eq:bidomain_model} and \eqref{Eq:monodomain_model}.

For the space discretization, we build two independent tetrahedral meshes on $\Omega_H$ and $\Omega_T$. Piecewise linear finite elements \cite{Ethier2008,QuarteroniValli2008} are used to solve both problems.  Additionally, we employ a piece-wise linear interpolation operator $\Pi$ at the heart-torso interface, given by
$$ u_T = \Pi u_e \quad \text{on }\Gamma \times (0,T],$$
which enables non-conforming spatial discretizations on the cardiac domain $\Omega_H$ and the torso domain $\Omega_T$ (see Fig. \ref{Fig:scheme_num} and \ref{Appendix} for a schematic representation and a more detailed description of the numerical scheme implemented, respectively.)

\begin{rmk}
	By employing linear interpolation on the heart-torso interface $\Gamma$, our approach accommodates both mesh-based non-conformity and geometrical non-conformity. Mesh-based non-conformity refers to the ability to have different discretizations on $\Gamma$ while maintaining the same shape and position on both $\Omega_H$ and $\Omega_T$. Geometrical non-conformity, on the other hand, allows for variations in the shape of the cardiac domain or the movement of the heart in space. In Section \ref{Sec:numerical_results}, we demonstrate the benefits of coarsening the torso mesh (mesh-based non-conformity) and explore the performance of our method in the presence of geometrical non-conformity.
\end{rmk}

\subsection{Computing clinical outputs}
\label{Subsec:computing_clinical_outputs}
\begin{figure}[!t]
	\centering
	\includegraphics[width=0.7\textheight]{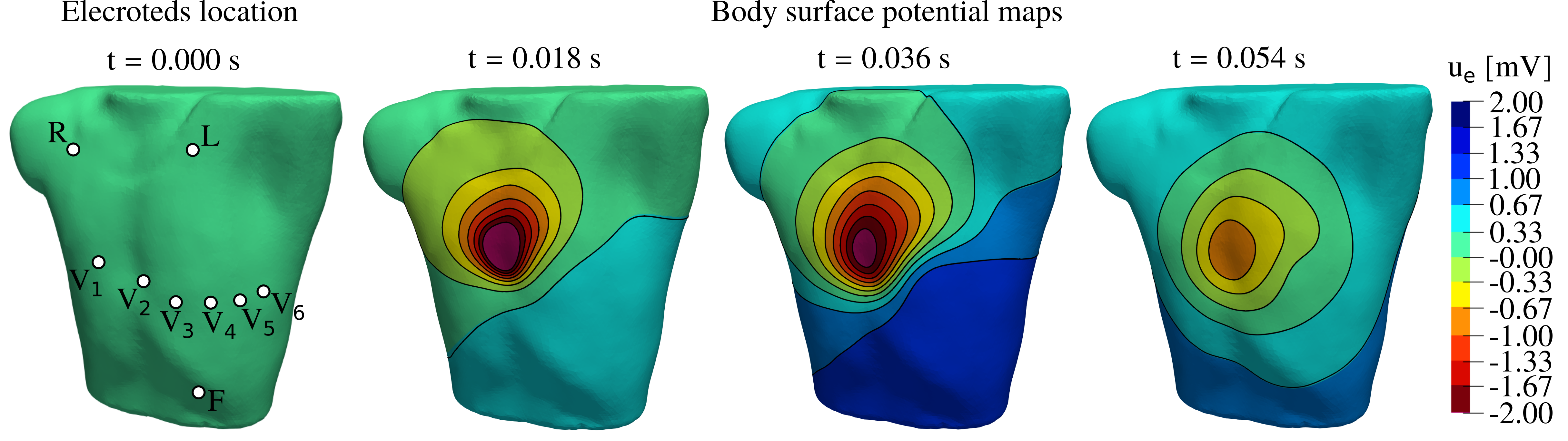}
	\caption{Electrodes location (left) and body surface potential maps on the torso domain $\Omega_T$, sample at time instants $t = 0.018, 0.036, 0.054 s$.}
	\label{Fig:ref_bspm}	
\end{figure}

The solution of the torso model \eqref{Eq:torso} provides the evolution in time of the extracellular potential $u_T$, spatially distributed on the whole torso tissue. Non-invasive clinical outputs as the BSPMs and the ECGs are, therefore, obtained by post-processing the outputs of the torso simulation. 
Body surface potential maps (BSPMs) can be computed by extracting the values of $u_T$ on the external surface $\Sigma^{\text{ext}}$ on fixed time instants (as shown in Fig. \ref{Fig:ref_bspm}). Electrocardiogram (ECG) signals, on the other hand, are time dependent signals obtained by observing the temporal evolution of $u_T$ at specific locations on the surface of the human body where electrodes are placed. Among the most frequently used ECG system, the numerical results of Section \ref{Sec:numerical_results} refer to the 12-lead ECG system \cite{Malmivuo1995,Halhuber1979}.

The standard 12-lead ECG is a system of 12 leads obtained by combining the values of $u_T$ recorded from 9 electrodes, named $R$, $L$, $F$, and $V_i$, $i = 1,\dots,6$ (see Fig. \ref{Fig:ref_bspm}). Defining by $\mathbf{x}_R$, $\mathbf{x}_L$, $\mathbf{x}_F$, $\mathbf{x}_{V_i}$, the spatial location of the electrodes, the 6 limb leads are computed as
\begin{gather*}
I = u_T(\mathbf{x}_L) - u_T(\mathbf{x}_R), \quad II = u_T(\mathbf{x}_F) - u_T(\mathbf{x}_R), \quad III = u_T(\mathbf{x}_F) - u_T(\mathbf{x}_L)\\ 
aVR = u_T(\mathbf{x}_R) - \frac{1}{2} (u_T(\mathbf{x}_L) + u_T(\mathbf{x}_F)) \quad aVL = u_T(\mathbf{x}_L) - \frac{1}{2}(u_T(\mathbf{x}_R) + u_T(\mathbf{x}_F)) \\ aVF = u_T(\mathbf{x}_F) - \frac{1}{2}(u_T(\mathbf{x}_L) + u_T(\mathbf{x}_R)),
\end{gather*}
whereas the 6 precordial (or chest) leads are defined as
\[ V_i = u_T(\mathbf{x}_{V_{i}}) -WCT, \quad i = 1,\dots,6,
\]
with $WCT$ denoting the Wilson central terminal signal, given by
$$WCT =  \frac{1}{3} \left[ u_T(\mathbf{x}_L) + u_T(\mathbf{x}_R) + u_T(\mathbf{x}_F)\right ].$$

\section{Numerical results}
\label{Sec:numerical_results}
The numerical simulations presented in this section are run on a basal biventricular geometry representing the cardiac domain $\Omega_H$, obtained by cutting with a plane the Zygote 3D human heart \cite{Zygote} (see Fig. \ref{Fig:EP_data}). This CAD model represents the heart of an average 21-year-old male. The torso domain $\Omega_T$ is instead the realistic RIUNET torso geometry \cite{Ferrer2015}(\url{http://hdl.handle.net/10251/55150}), obtained from the online repository at the Center for Integrative Biomedical Computing \cite{CIBC} of the University of Utah. The original atria, ventricles and blood pools have been removed and replaced with the Zygote geometry. Processing and meshing of $\Omega_H$ and $\Omega_T$ are carried out by means of the software library VMTK \cite{Antiga2008} and the open access Paraview software \cite{paraview}. Solvers for the coupled model were implemented in the C++ library \texttt{life$^x$} \cite{africa2022lifex,africa2023fibers,lifex_2} based on the finite element core \texttt{deal.II} \cite{dealII93,arndt2021deal,dealii}. All the simulation are run at the CINECA high-performance computing center (Italy), employing 192 parallel processes on the GALILEO100 supercomputer \footnote{see \url{https://wiki.u-gov.it/confluence/display/SCAIUS/UG3.3\%3A\%2BGALILEO100\%2BUserGuide} for the technical specifications.}.

We introduce a tetrahedral mesh on both cardiac $\Omega_H$ and torso $\Omega_T$ domains, initially ensuring conforming discretization at the heart-torso interface $\Gamma$. 
The mesh constructed on the basal biventricular region, $\Omega_H$, has an average resolution of $1.05 mm$, resulting in a total of 3130352 vertices. On the torso domain we employ a coarser discretization by creating a gradient mesh between the already discretized surfaces, $\Gamma$ and $\Sigma^{\text{ext}}$. To ensure interface conformity, we extract $\Gamma$ from $\Omega_H$, while we build a surface mesh with an average size of $5.00 mm$ on $\Sigma^{\text{ext}}$. Consequently, the average mesh size in $\Omega_T$ is reduced to $2.53 mm$, with 1471944 vertices. Finally, we use second-order finite elements to increase the number of degrees of freedom in both $\Omega_H$ and $\Omega_T$, while simultaneously halving the mesh resolution.
We referred to these initial spatial discretizations as \emph{reference} cardiac and torso meshes.

We simulate 4 heart-beats with $T = 3.2s$, and consider only the last heartbeat for our analysis (starting from $t_0 = 2.4s$), to reduce the effect of the initial conditions of the ionic model. The simulation time step is set to $\Delta t_H = 50 \mu s$ for the EP simulation. A larger time step of $\Delta t_T = 20 \Delta t_H = 1 ms$ is instead employed to solve the elliptic model in the torso.

The 12-lead ECG system, as well as the body surface potential maps, are finally computed post-processing of each simulation, according to subsection \ref{Subsec:computing_clinical_outputs}. Since the cardiac geometry only represents the two ventricles, the P wave cannot be computed. Moreover, we do not model cell heterogeneities in the cardiac tissue, which is responsible for the physiological representation of the T wave \cite{Keller2012,XUE2010553,Weiss2007}. The ECG analysis is thus restricted to the QRS complex (\emph{e.g.} see Fig. \ref{Fig.variability_methods}). 

\subsection{Model comparison}
\begin{figure}[!t]
	\centering
	\includegraphics[width=0.65\textheight]{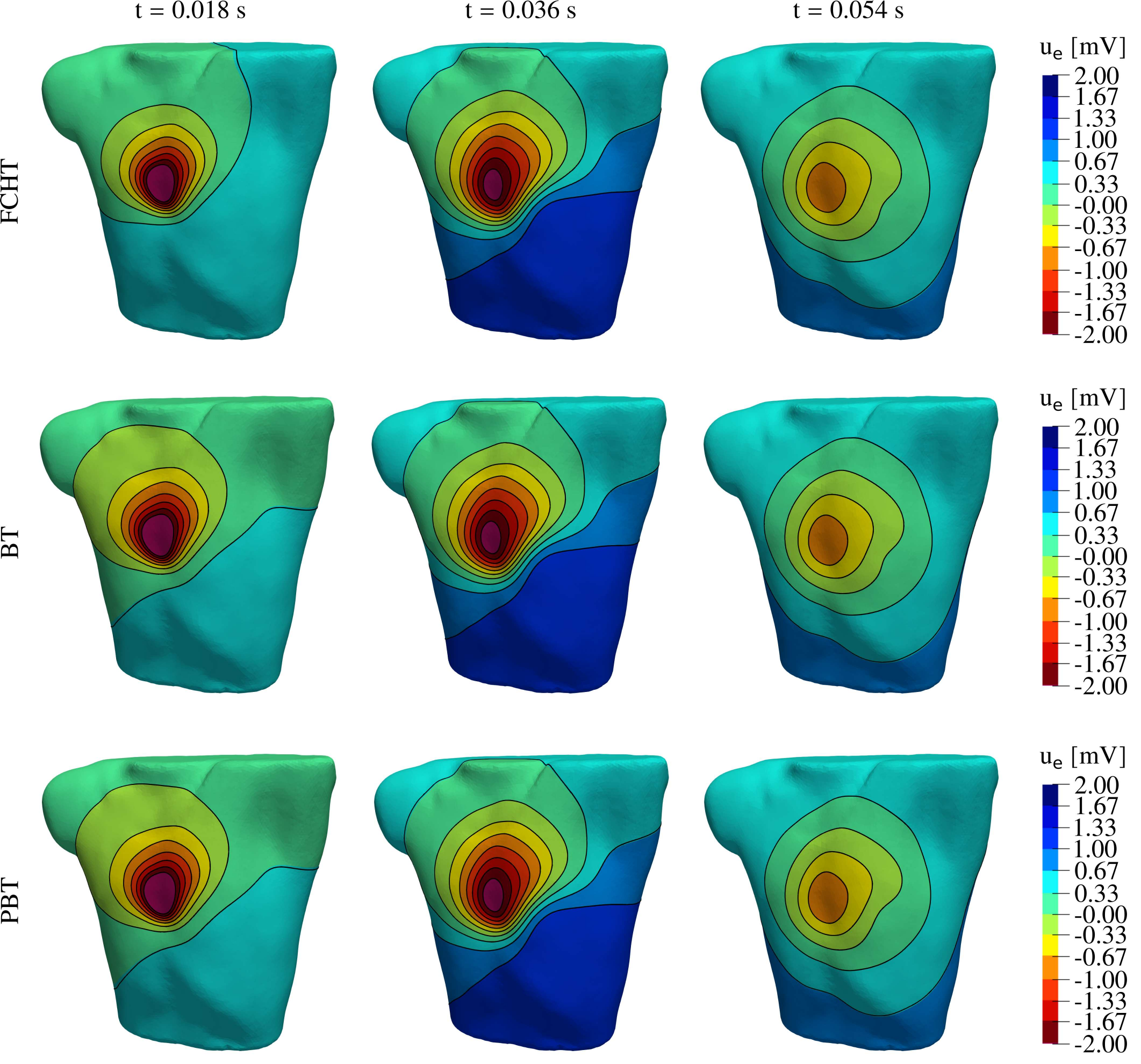}
	\caption{Body surface potential maps computed with the fully-coupled heart-torso (FCHT) model, the bidomain-torso (BT) model, and pseudo-bidomain-torso (PBT) model, sample at time instants $t = 0.018, 0.036, 0.054 s$.}
	\label{Fig.torso_BSPM_methods}
\end{figure}

\begin{figure}[!t]
	\centering
	\includegraphics[width=0.72\textheight]{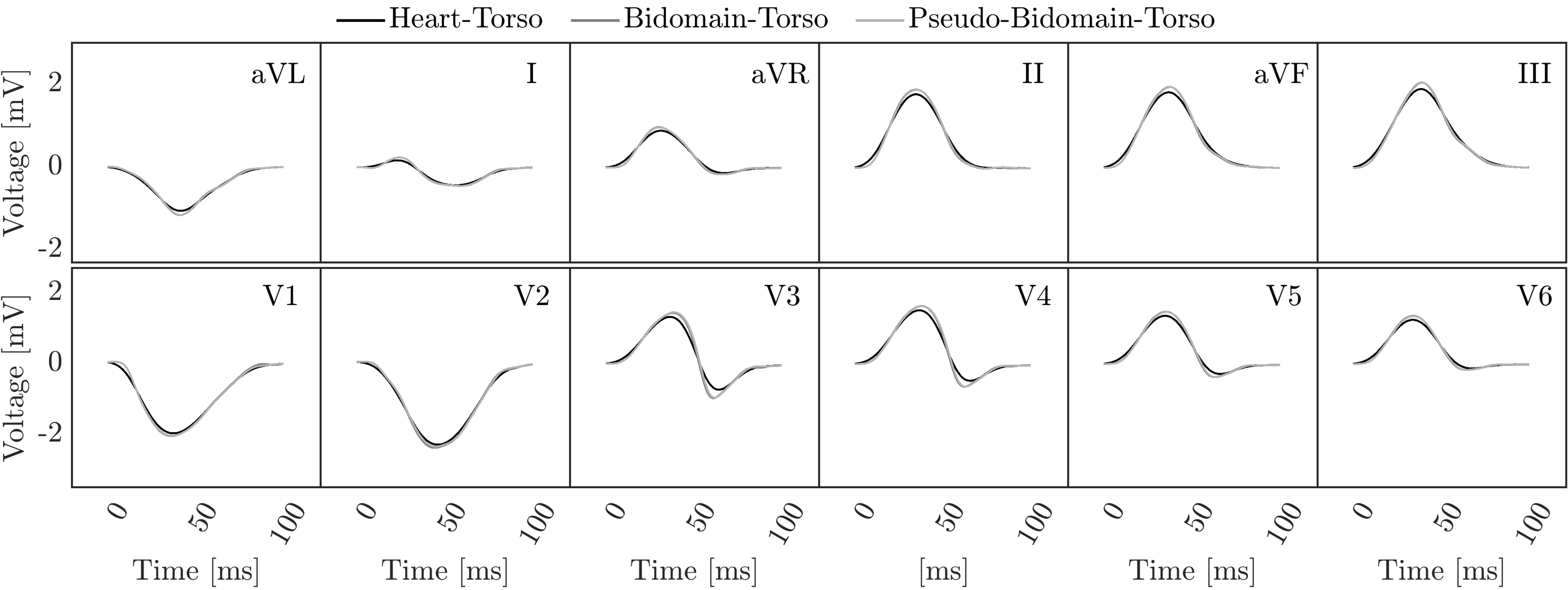}
	\caption{QRS complex computed with the FCHT model, the bidomain-torso model, and the  pseudo-bidomain-torso model.}
	\label{Fig.variability_methods}
\end{figure}

In order to assess the reliability and efficiency of the proposed framework, we compare the FCHT model, the bidomain-torso model, and the pseudo-bidomain torso model. We analyze the difference between the models solutions both qualitatively in Fig. \ref{Fig.torso_BSPM_methods} and \ref{Fig.variability_methods}, and quantitatively in Table \eqref{Tab:variability_models}. Body surface potential maps are quantitatively compared by computing a relative $\ell_2$ norm over the body surface, averaged over time.
At each lead, differences between signals $\phi_1$ and $\phi_1$ are determined with the relative root-mean-square error (rmse)
\begin{equation}
\label{Eq:RMSE}
\displaystyle 
\text{rmse} = \frac{\sqrt{\frac{1}{N_T}\sum_{i=1}^{N_T} \left[\phi_1(i) - \phi_2(i)\right]^2}}{\sqrt{\frac{1}{N_T}\sum_{i=1}^{N_T} \phi_1(i)^2}},
\end{equation}
accounting for variations of both signal morphology and relative amplitude, where $N_T$ is the number of time step in which the torso model \eqref{Eq:torso} is solved. The overall variance on the 12-lead system is then accounted for by computing the mean of the rmse values. However, to study the ECG shape variation independently of the amplitude scaling, we also compute the mean of the linear correlation coefficient
\begin{equation}
\label{Eq:cc}
\text{CC} = \frac{1}{s_1 s_2} \sum_{i=1}^{N_T}\left[\phi_1(i) - \overline{\phi_1}\right]\left[\phi_2(i) - \overline{\phi_2}\right]
\end{equation}
of each lead. Here $s_{1,2}$ are the standard deviations, and $\overline{\phi}_{1}$ and $\overline{\phi}_{2}$ are the arithmetic mean values of $\phi_{1}$ and $\phi_2$, respectively.

The variations in body surface potential maps resulting from different model assumptions are most noticeable in the initial time instants of the simulation (see Fig. \ref{Fig.torso_BSPM_methods}), but diminish over time. This can be attributed to the fact that depolarization begins in the right ventricle, which is the closest part of the cardiac domain to the anterior part of the body, therefore making the BSPM more susceptible to variations in EP solutions. Nevertheless, quantitatively, the average differences between the maps remain below 5\% (see Table \ref{Tab:variability_models}), regardless of the model being considered.

Nearly identical ECGs are obtained using the bidomain-torso and pseudo-bidomain-torso models, while minor signal variations are observed when comparing the ECGs generated by the one-way coupled model with those of the FCHT model. It is noteworthy that both the shape and magnitude of the signals remain largely unaffected by the mathematical model used (refer to Figure \ref{Fig.variability_methods}). Indeed, only slight variations are shown by the signal amplitude in leads II, aVF, and III, as well as a more negative S wave in leads V$_3$ and V$_4$. However, these variations are minimal and can be qualitatively disregarded when simulating healthy conditions.

The qualitative findings are coherent with the quantitative metrics presented in Table \ref{Tab:variability_models}. The average root mean square error (RMSE) is consistently below 10\%, regardless of the specific one-way electrophysiology-torso model utilized. Additionally, the correlation coefficients show that significant variations in the ECGs are primarily attributed to amplitude scaling rather than signal morphology.
\begin{table}[!t]
	\centering 
	\begin{tabular}{cccccccc}
		\toprule
		&$~$BSMP$_{\ell_2,mean}$ &rmse$_{\text{mean}}$ &CC$_{\text{mean}}$ &\% of costs $~$ &Speed up &\multicolumn{2}{c}{Partial \% of costs}\\
		\hline &&&&& \\ [-2ex]
		FCHT &&&&100\% &&&\\
		\hline &&& \\ [-2ex]
		\multirow{3}{*}{BT} &\multirow{3}{*}{3.55\%} &\multirow{3}{*}{9.98\%} &\multirow{3}{*}{0.996} &\multirow{3}{*}{59\%} &\multirow{3}{*}{1.7} &Bidomain &91.3\% \\ &&&&&&Torso &4.5\% \\ &&&&&&Interface interpolation &4.2\% \\
		\hline &&&&&& \\ [-2ex]
		\multirow{3}{*}{PBT} &\multirow{3}{*}{5.28\%} &\multirow{3}{*}{9.65\%} &\multirow{3}{*}{0.997}  &\multirow{3}{*}{48\%} &\multirow{3}{*}{2.1} &Pseudo-bidomain &91.6\% \\ &&&&&&Torso &4.0\% \\ &&&&&&Interface interp&4.4\% \\
		\bottomrule
	\end{tabular}
	\caption{Second column: relative $\ell_2$ error between the body surface potential maps (BSMP) computed with the FCHT, BT and PBT models, averaged over time. Third column: rmse errors computed between the FCHT model and either the BT model or the PBT model, averaged over the 12 leads. Fourth column: correlation coefficients between the ECG signals computed with either the BT and the FCHT models, and the PBT and FCHT models. Fifth column: comparison of the computational costs associated with solving the FCHT, the BT, and the PBT models. We present the percentage ratio between each segregated-staggered scheme and the fully-coupled model, which serves as the reference model. We indicate the percentage ratio between each segregated-staggered scheme with the fully-coupled model, which is here considered as the reference model.
	Sixth column: Computational cost speedup achieved with the segregated-staggered algorithms.
	Seventh column: Partial computational costs expressed as percentages. We consider three main phases: the solution of the electrophysiological model, the solution of the torso model, and the pre-processing of the Dirichlet interface conditions.}
	\label{Tab:variability_models}
\end{table}	

The efficiency of the proposed computational framework is evaluated by analyzing the associated costs, as presented in Table \ref{Tab:variability_models}. Our results demonstrate that the one-way models offer significantly improved cost-effectiveness compared to the FCHT model, with a reduction in simulation time exceeding 50\% when employing the pseudo-bidomain-torso model. In addition, we partition the computational cost of the simulation into three fractions, each representing the time required to solve the core components of the coupled model: the electrophysiology model, the torso model, and the signal interpolation at the interface. This breakdown allows us to examine the impact of incorporating signal interpolation more closely. Our analysis reveals that the cost of interpolation is comparable to that of the torso solution, accounting for only 5\% of the overall simulations. Thus, the additional cost associated with signal interpolation is minimal in comparison to the total simulation cost, which is primarily driven by the expense of solving the electrophysiology model. By leveraging the implemented interpolation method, we will show that the simulation costs can be further reduced by decreasing the mesh complexity within the torso domain, meanwhile preserving result accuracy. The efficacy of this approach will be demonstrated in the following subsection.

\subsection{Mesh-based interface non-conformity}
In this section we aim at assessing the major potential of our segregated, staggered and interface interpolation based computational framework, by demonstrating its ability to preserve the results accuracy while significantly reducing the computational costs associated with solving the torso model, that is to investigate the effect mesh-based non-conformity between the heart and torso domains. This is accomplished by gradually reducing the mesh resolution in the torso domain, and by projecting a prescribed electrophysiology solution onto the resulting torso meshes. Throughout this analysis, the prescribed electrophysiology solution is kept consistent for each new torso simulation, as it is obtained by solving the electrophysiology model with identical parameter settings and mesh resolution in the cardiac domain. This analysis stands as a noteworthy contribution in its own right, as it explores the impact of torso discretization on the ECG and body potential maps, highlighting the negligible influence of such discretization on the overall model results.

\begin{figure}[!t]
	\centering
	\includegraphics[width=0.75\textheight]{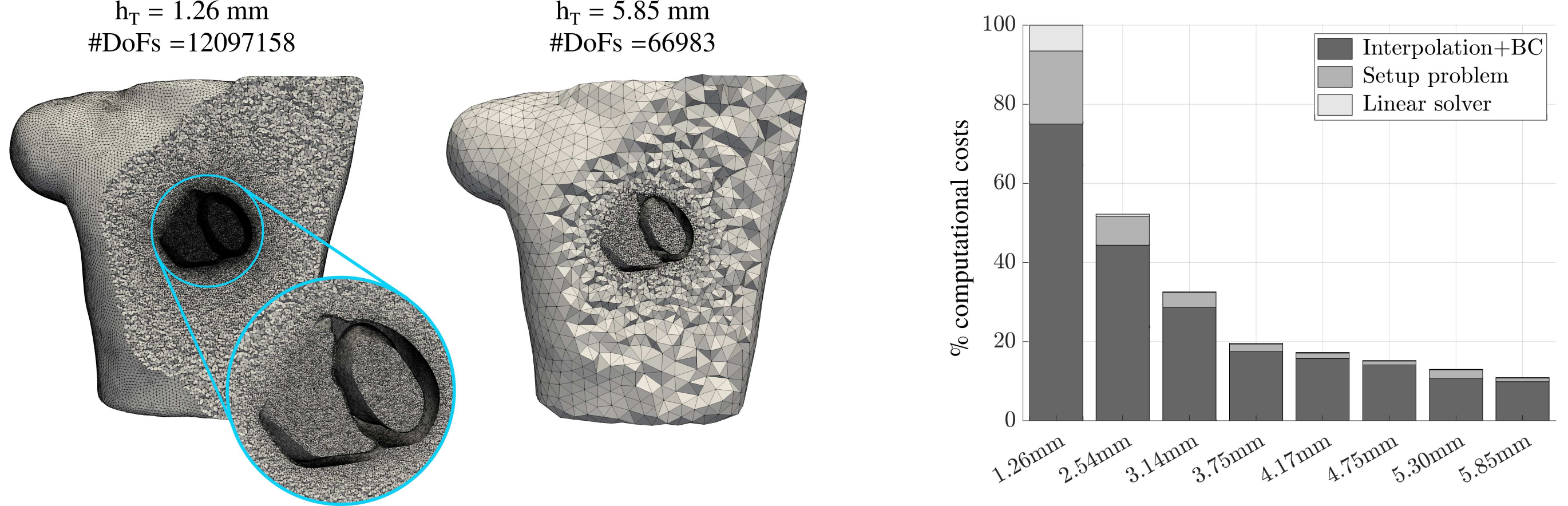} 
	\caption{Left: Comparison of finer and coarser discretizations of the torso domain.
		Right: Computational costs associated with solving the diffusive model in the torso using different discretizations of the torso geometry. Each bar represents the costs of a single simulation, shown as a percentage ratio relative to the reference simulation ($h = 1.26mm$). The bars also indicate the distribution of costs across three main phases: system setup (setup problem), interpolation and application of the Dirichlet interface conditions (interpolation + BC), and solution of the linear system (linear solver).}
	\label{Fig.info_mesh_sizes}
\end{figure}

Based on the findings presented in the previous subsection, we employ the pseudo-bidomain-torso model for the current study. Seven different discretizations are built on the torso domain. As for the reference framework, a gradient mesh is generated within the torso domain, beginning with prescribed discretizations on the interface $\Gamma$ and the external surface $\Sigma^{\text{ext}}$, spanning from an average mesh size of 1.26mm, corresponding to 12097158 DoFs, to an average mesh diameter of 5.85mm with 66983 DoFs (see Fig. \ref{Fig.info_mesh_sizes}, right). The mesh sizes of both surfaces, as well as the average size of the final gradient meshes, and total number of DoFs for each torso mesh are documented in Table \ref{Tab:variability_meshes}.

\begin{table}[!t]
	\centering 
	\begin{tabular}{cccc|ccc|ccc}
		\toprule
		$h_{T,\Gamma}$ &$h_{T,\Sigma^{\text{ext}}}$ &$h_{T}$ &\#DoFs &rmse$_{\text{mean}}$ &rmse$_{V_1}$ &rmse$_{V_2}$ &CC$_{\text{mean}}$ &CC$_{V_1}$ &CC$_{V_2}$ \\
		\hline 
		&&&&&&&&& \\ [-2ex]
		0.50 mm &2.50 mm &1.26 mm &12097158&&&&&&\\
		1.00 mm &5.00 mm &2.54 mm &1471944 &1.7\% &3.6\% &0.5\% &0.9999 &0.9999 &1.000\\
		1.25 mm &7.50 mm &3.14 mm &610733 &3.8\% &13.8\% &7.0\% &0.9994 &0.9974 &0.9970\\
		1.50 mm &10.00 mm &3.75 mm &314643 &3.1\% &5.7\% &9.7\% &0.9994 &1.000&0.9938\\
		1.75 mm &12.50 mm &4.17 mm &204199 &6.2\%&23.0\% &11.5\% &0.9987 &0.9977&0.9941\\
		2.00 mm &15.00 mm &4.75 mm &131655 &6.6\% &9.4\% &12.5\% &0.9982 &0.9984&0.9895\\
		2.25 mm &17.50 mm &5.30 mm &91636 &7.1\% &14.5\% &23.6\% & 0.9974 &0.9932&0.9839\\
		2.50 mm &20.00 mm &5.85 mm &66983 &6.9\% &22.7\% &13.9\% & 0.9974 &0.9860&0.9909\\
		\bottomrule
	\end{tabular}
	\caption{Average mesh size of the heart-torso interface $\Gamma$, of the external surface of the torso $\Sigma^{\text{ext}}$, and of the torso mesh. We moreover report the rmse and CC computed between the reference signal, \emph{i.e.} the ECG obtained by employing the finest mesh in the torso (row 1), and the signal test, those obtained by coarsening the torso discretization, average over all the 12 leads, and for leads $V_1$ and $V_2$.}
	\label{Tab:variability_meshes}
\end{table}
  
\begin{figure}[!t]
	\centering
	\includegraphics[width=0.72\textheight]{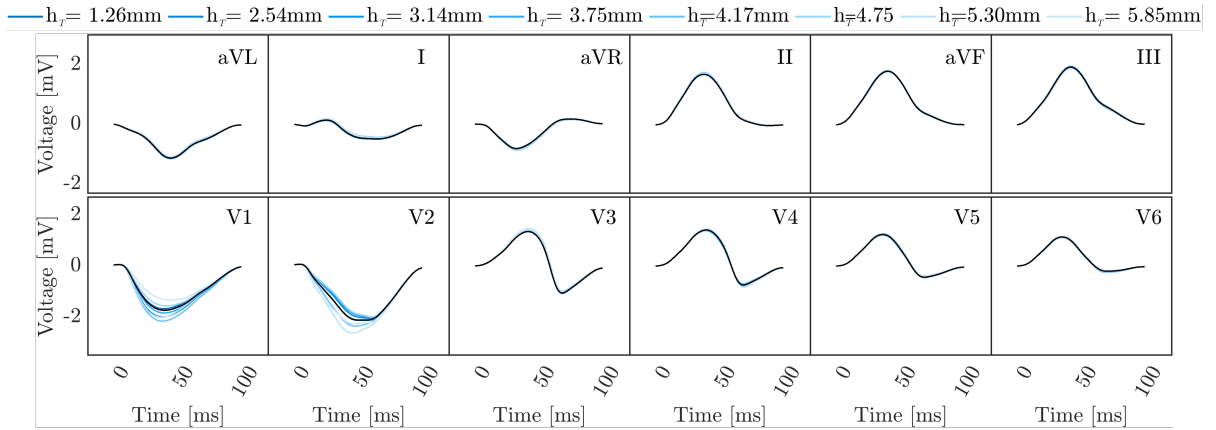}
	\caption{QRS complex computed employing different discretizations in the torso geometry. The legend reports the average mesh sizes.}
	\label{Fig.variability_mesh_size}
\end{figure}

\begin{figure}[!t]
	\centering
	\includegraphics[width=0.65\textheight]{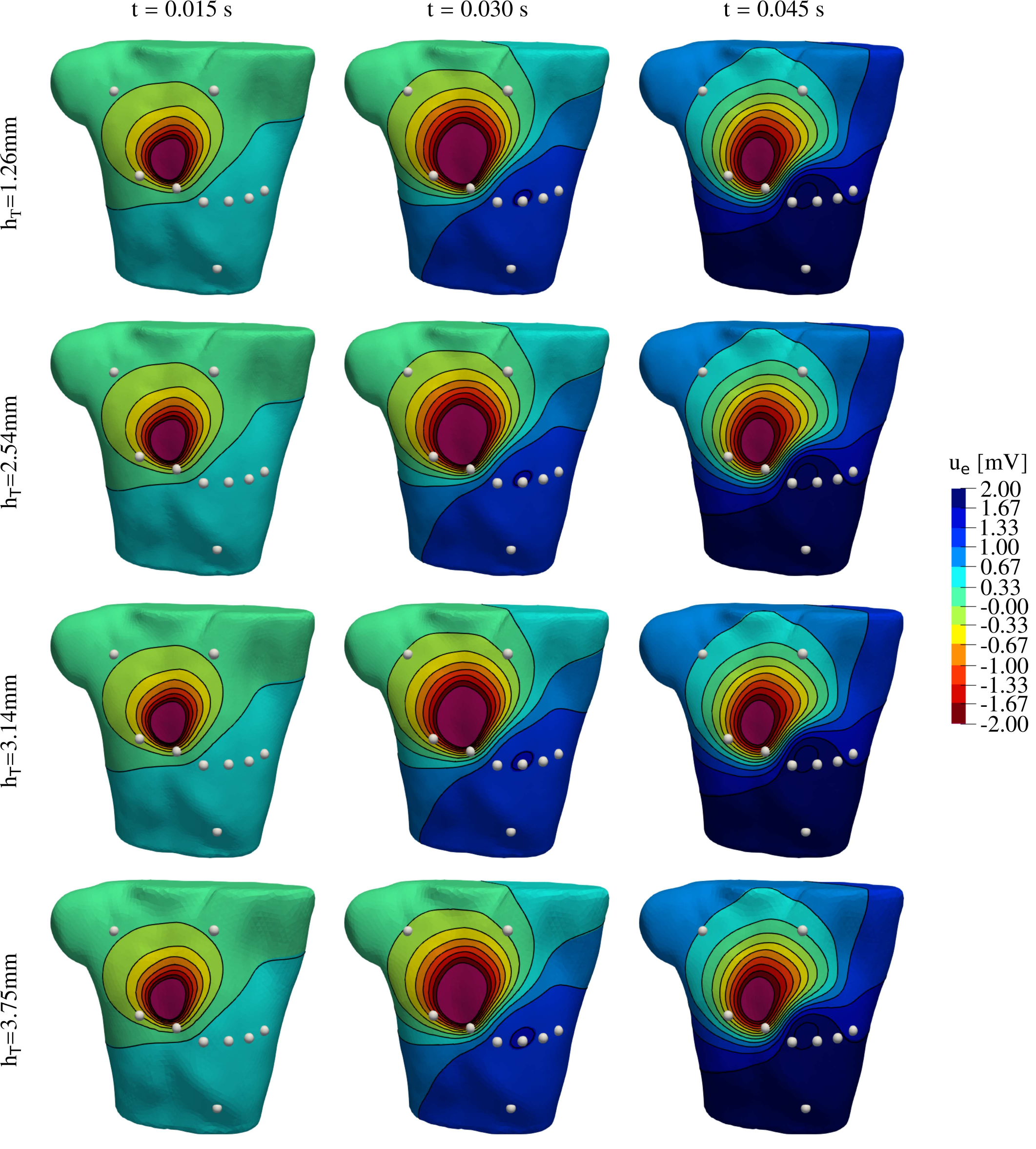}
	\caption{Body surface potential maps computed on different torso domain discretizations, with average mesh size between $1.26mm$ and $3.75mm$, sampled at time instants $t = 0.015, 0.030, 0.045s$.}
	\label{Fig.Torso_mesh_size_1}
\end{figure}

\begin{figure}[!t]
	\centering
	\includegraphics[width=0.65\textheight]{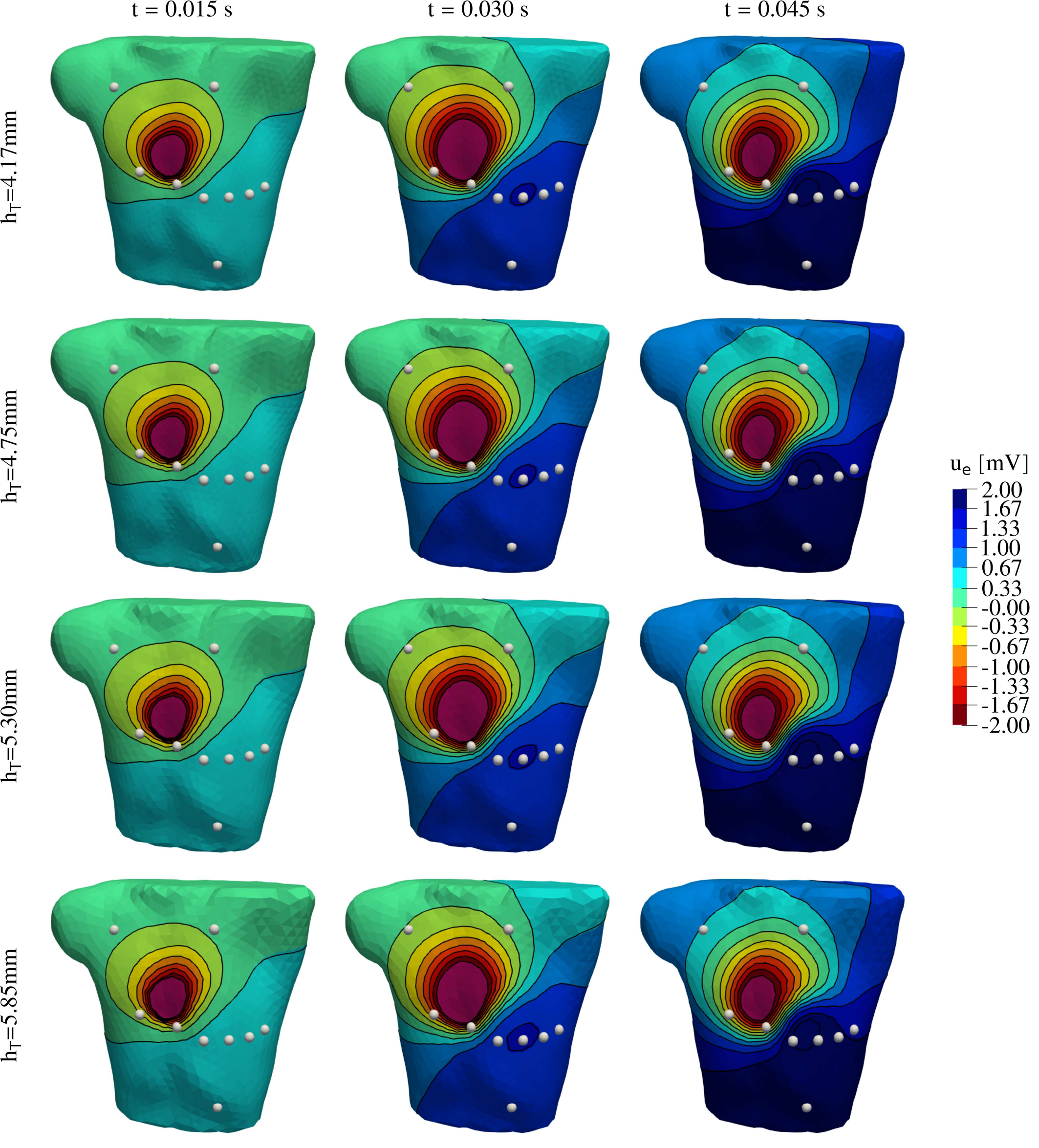}
	\caption{Body surface potential maps computed on different torso domain discretizations, with average mesh size between $4.17mm$ and $5.85mm$, sampled at time instants $t = 0.015, 0.030, 0.045s$..}
	\label{Fig.Torso_mesh_size_2}
\end{figure}

Body surface potential maps remain unaffected by changes in the discretization of the torso domain (see Fig. \ref{Fig.Torso_mesh_size_1}-\ref{Fig.Torso_mesh_size_2}). The ECG traces exhibit similar patterns across all leads, with slight amplitude variations observed in leads V1 and V2, as shown in Fig. \ref{Fig.variability_mesh_size}. However, such variations can be expected, since these leads are in close proximity to the electric sources represented by the heart, making them more susceptible to the DoFs spatial distribution within the torso domain. In quantitative terms, the average rmse and CC reported in Table \ref{Tab:variability_meshes} validate these findings. The ECGs exhibit a gradual variation, primarily in amplitude, as the mesh size becomes coarser, resulting in an average RMSE of 7\% in the coarser discretization. Notably, higher values are found for leads $V_1$ and $V_2$, with trace variations of up to 23\% compared to the reference finer torso mesh. However, despite these amplitude variations, the correlation coefficients indicate that the changes in morphology remain minimal, even for leads V1 and V2.

Through this experiment, we demonstrate that BSPMs and ECGs remain nearly identical regardless of the discretization of the torso domain. Consequently, the mesh resolution in the torso domain can be significantly coarsened to minimize the simulation burden. In comparing computational consumption, we focus solely on the solution of the torso problem over time. The total computational costs are divided into three components: problem setup, problem solving, and signal interpolation at the domain interfaces. Using the simulation with the finest mesh in the torso as a reference, we observe a rapid decrease in computational costs, with a reduction of up to 90\% in total simulation costs when employing the coarsest discretization, being the interpolation scheme to be accounted for the majority of the computational costs.

\subsection{Geometrical interface non-conformity}
The final numerical experiment aims at assessing the flexibility of the proposed numerical framework concerning geometrical interface non-conformity. Specifically, it investigates the reliability of the BSPMs and ECGs results when combining a cardiac and torso geometry with non-matching interfaces in terms of spatial alignment. Geometrical non-conformity can arise in model testing scenarios or when model output features related to geometrical properties of the cardiac domain are under investigation, where employing or computing a patient-specific geometry may be prohibitively expensive. In such cases, the ability to freely modify the cardiac geometry while preserving the torso geometry, without the need for domain modifications or remeshing, can significantly reduce computational time. This study is then carried out by applying rigid transformations to the cardiac domain while keeping the torso domain in the reference configuration described earlier. This approach effectively generates a set of heart-torso domains with geometric non-conformity.

\begin{figure}[!t]
	\centering
	\includegraphics[width=0.7\textheight]{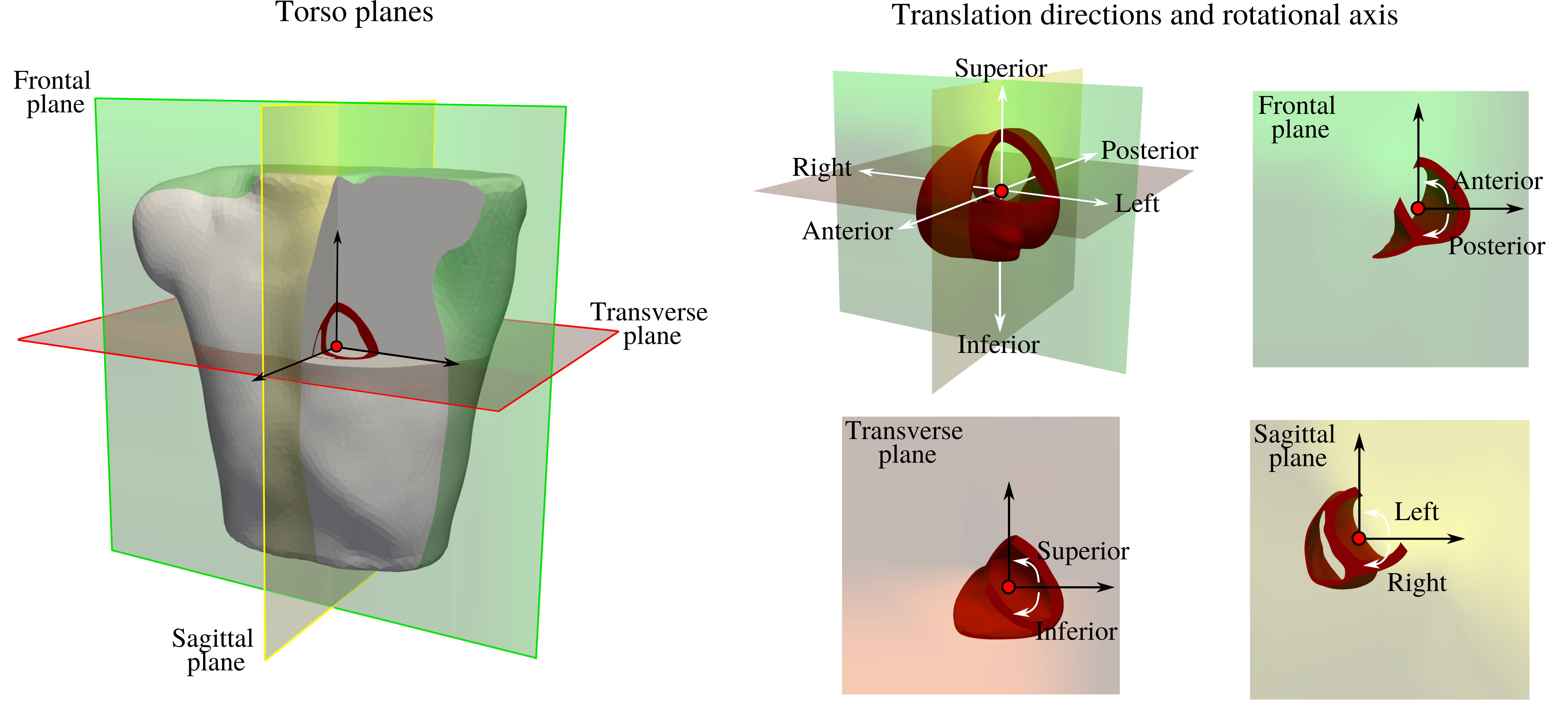} 
	\caption{Left: sagittal, frontal and transverse plane on the torso domain. Right: graphical representations of the six translation directions and six corresponding rotational axes along which the cardiac domain is displaced.}
	\label{Fig.Torso_planes}
\end{figure}

All simulations are performed with the pseudo-bidomain-torso model. We prescribed mesh with average mesh size $1.5mmm$ and $10.0mm$ on the interface $\Gamma$, and the external surface of the torso $\Sigma^{\text{ext}}$, respectively. An orthogonal frame is established, centered at the barycenter of the heart, with axes defined by the intersections of the torso's sagittal, frontal, and transverse planes. By dividing each axis into two segments, we define six translation directions and six corresponding rotational axes along which the cardiac domain is displaced. Specifically, the translation directions include left and right (based on the intersection of the transverse and frontal planes), superior and inferior (based on the intersection of the frontal and sagittal planes), posterior and anterior (based on the intersection of the sagittal and transverse planes), see Fig. \ref{Fig.Torso_planes}). For each translation direction, we specify three translation measures, \emph{i.e.} 3cm, 6cm and 9cm, and for each rotational axis, we prescribe three rotation angles, \emph{i.e.} of 3, 6 and 9 degrees, resulting in a total of 36 heart-torso geometries. These geometries are further divided into two groups: one representing small transformations, such as those resulting from the respiratory motions of the heart \cite{shechter2004respiratory}, and another representing large transformations, intended to encompass the range of heart positions within the torso across a human population \cite{Odille2017}. Through experiments on the first group, our objective is to demonstrate the model's capability to accurately reproduce reliable results, both quantitatively and qualitatively, when considering the geometrical variability of the same patients in the simulations. In the second group of experiments, our aim is to showcase that our framework can qualitatively reproduce various clinically relevant outputs. In the first case, we present BSPMs and ECGs, along with quantitative indices such as rmse and CC between the signals computed with conforming and non-conforming geometries. In the latter case, we focus on a visual comparison of the ECG traces.

\begin{figure}[!t]
	\centering
	\includegraphics[width=0.65\textheight]{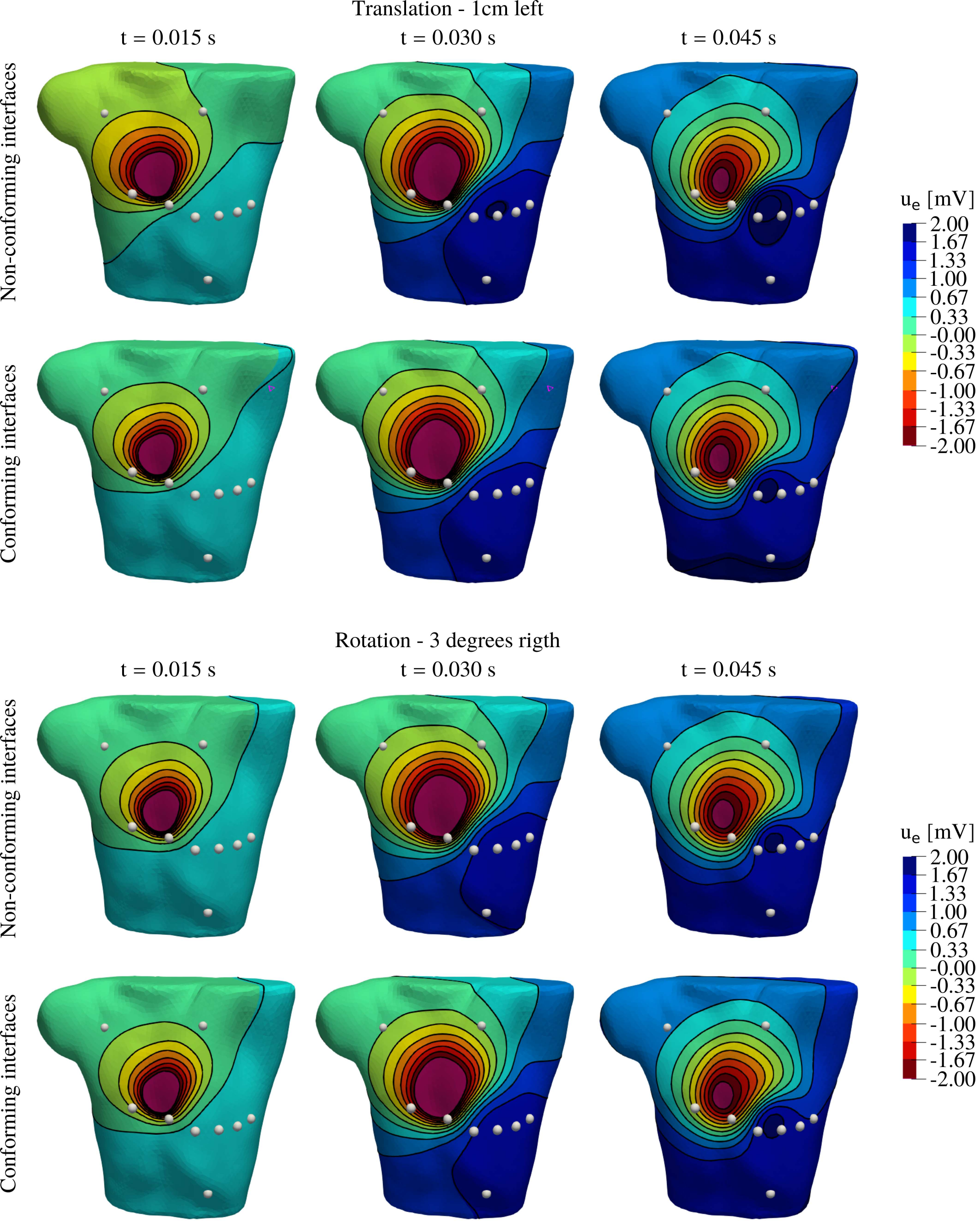}
	\caption{Body surface potential maps obtained by displacing the cardiac domain of $1cm$ on the left (top) and or 3 degree on the right (bottom). For each set of cardiac geometries, we compare the BSPMs on the reference conforming heart-torso geometry, and on geometrical non-conforming geometry sampled at time instants $t = 0.015, 0.030, 0.045s$.}
	\label{Fig.Torso_small_movements}
\end{figure}

The proposed numerical approach yieldss good performances in computing body surface potential maps, even when dealing with geometric non-conformity resulting from small transformations of the cardiac domain. It effectively preserves the key features of solution propagation, particularly in cases involving rigid rotations of the cardiac domain (refer to Figure \ref{Fig.Torso_small_movements} for a comparison of BSPMs in one rotation and translation scenario). 

\begin{figure}[!t]
	\centering
	\includegraphics[width=0.7\textheight]{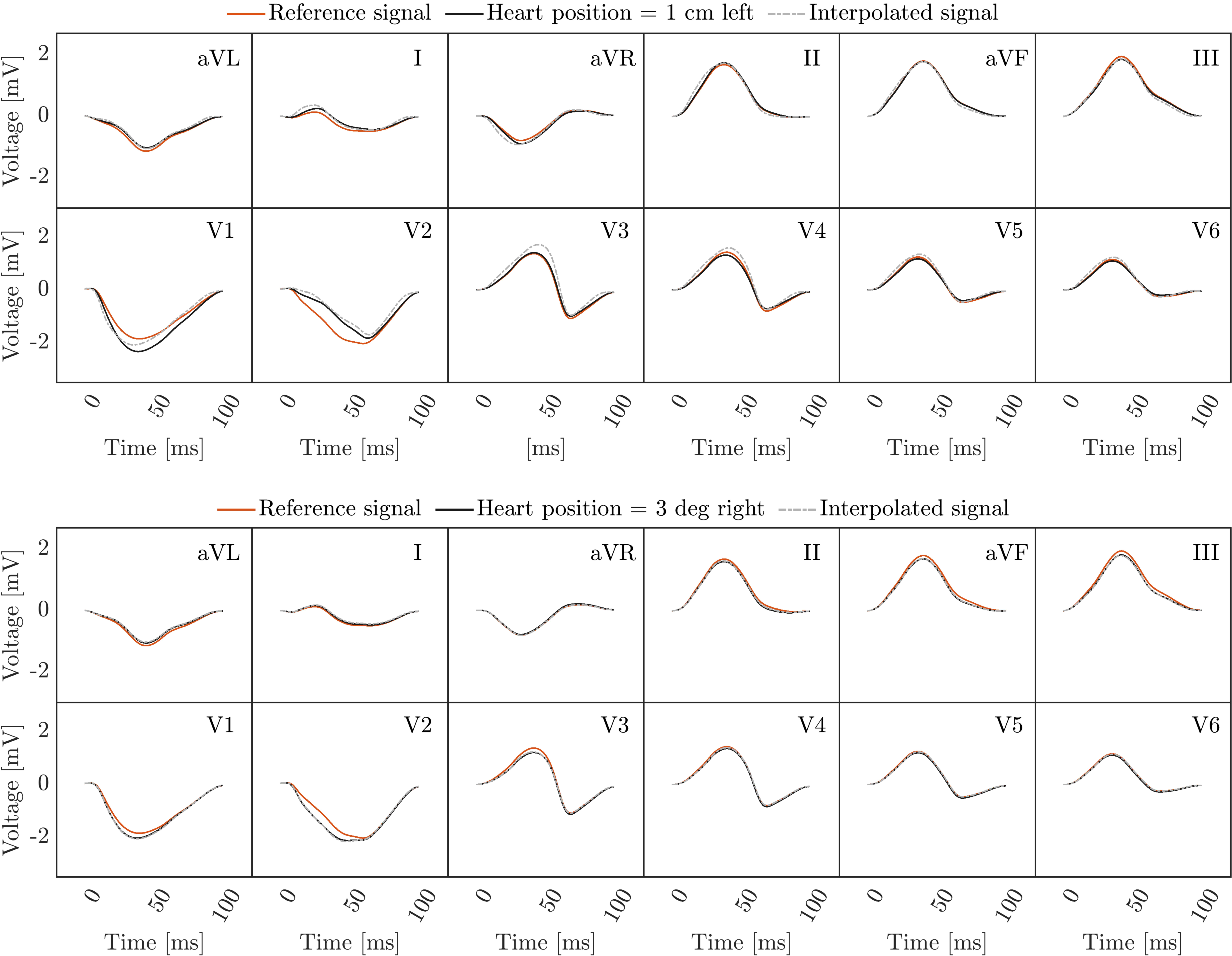} 
	\caption{Examples of QRS complex computed with reference heart-torso geometry (reference signal - red), and displacing the cardiac geometry but employing conforming heart-torso domains (heart position - black), and non-conforming domains (interpolated signal - grey). In this picture we analyze the effect of using the same torso geometry when subjecting the cardiac domain to small geometrical transformations \emph{e.g.}, translating the heart of 1 cm left (top), or rotating the heart of 3 degree around the left to right torso axis (bottom).}
	\label{Fig.examples_interpolation_error_heart_breath}
\end{figure}

\begin{table}[!t]
	\centering 
	\begin{tabular}{l|ccc|ccc}
		\toprule
		Transformation  &rmse$_{\text{mean}}$ &rmse$_{V_1}$ &rmse$_{V_2}$ &CC$_{\text{mean}}$ &CC$_{V_1}$ &CC$_{V_2}$ \\
		\hline 
		&&&&&& \\ [-2ex]
		Translation - 1cm right &17.0\% &12.9\% &33.1\% &0.9847 &0.9945 &0.9439\\
		Translation - 1cm left &15.4\% &15.2\% &29.3\% &0.9929 &0.9946 &0.9779\\
		Translation - 1cm superior &11.2\% &6.7\% &18.0\% &0.9955 &0.9978 &0.9982\\
		Translation - 1cm inferior &12.0\%&5.2\% &27.1\% &0.9956 &0.9990 &0.9769\\
		Translation - 1cm anterior &22.4\% &12.9\% &20.1\% &0.9590 &0.9971 &0.9967\\
		Translation - 1cm posterior &21.4\% &12.5\% &24.7\% & 0.9940 &0.9975 &0.9243\\
		Rotation - 3deg right &3.6\% &1.3\% &3.9\% &0.9994 &1.000 &0.9993\\
		Rotation - 3deg left &3.8\% &2.0\% &4.7\% &0.9994 &1.000 &0.9998\\
		Rotation - 3deg superior &4.4\% &1.2\% &5.5\% &0.9992 &0.9999 &0.9987\\
		Rotation - 3deg inferior &3.9\%&1.4\% &6.2\% &0.9992 &0.9989 &0.9798\\
		Rotation - 3deg anterior &2.6\% &1.5\% &2.8\% &0.9997 &0.9999 &0.9996\\
		Rotation - 3deg posterior &2.9\% &1.9\% &2.7\% & 0.9997 &0.9998 &0.9996\\
		\bottomrule
	\end{tabular}
	\caption{Quantitative indices comparing ECGs obtained from the pseudo-bidomain-torso model on geometrically non-conforming heart-torso geometries and conforming ones, obtained by applying small rigid transformations to the cardiac domain. We considered the rmse and the CC averaged over the 12 leads, and for leads $V_1$ and $V_2$.}
	\label{Tab:variability_geometry}
\end{table}

Similar results are also observed when analyzing the ECG traces (see Fig. \ref{Fig.examples_interpolation_error_heart_breath}). When non-conforming interfaces are employed, the ECGs are accurately reproduced, particularly when the heart undergoes rotation, resulting in slight variations, primarily in amplitude, in the precordial leads, such as V1 and V2. Table \ref{Tab:variability_geometry} provides quantitative indices for each of the 12 small transformations. For cardiac translations, the root mean square error is approximately one order of magnitude higher than in the other experiments discussed in this paper. However, the CCs indicate that, once again, major variations are primarily attributed to amplitude differences rather than changes in waveform shape. 

\begin{figure}[!t]
	\centering
	\includegraphics[width=0.65\textheight]{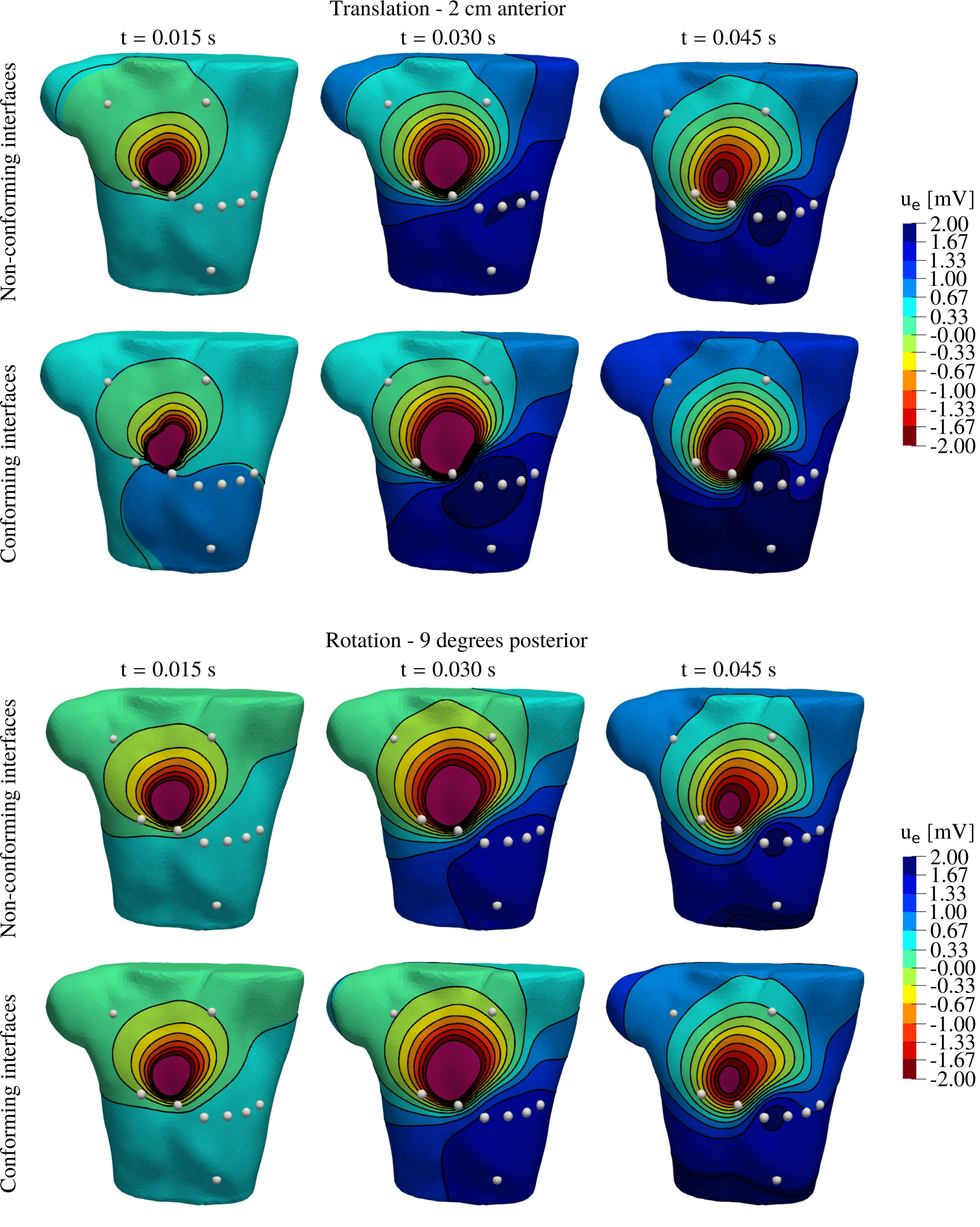}
	\caption{Body surface potential maps obtained by displacing the cardiac domain of $2cm$ on the directed toward the anterior of the torso domain (top) and or 3 degree directed to the posterior of the torso domain (bottom). For each set of cardiac geometries, we compare the BSPMs on the reference conforming heart-torso geometry, and on geometrical non-conforming geometry,sampled at time instants $t = 0.015, 0.030, 0.045s$.}
	\label{Fig.Torso_big_movements}
\end{figure}

While the proposed approach effectively reproduces desired outputs when there is slight difference between the heart and torso domains, increasing the degree of geometric non-conformity significantly challenges the accurate reproduction of model outputs. However, we are still able to faithfully reproduce the BSPMs when considering larger rotations of the cardiac domains, as shown in Fig. \ref{Fig.Torso_big_movements}. Instead, remarkable differences in the surface signal emerge when considering translations, indicating the need for a more careful approach in such cases (see \ref{Fig.Torso_big_movements}).

\begin{figure}[!t]
	\centering
	\includegraphics[width=0.6\textheight]{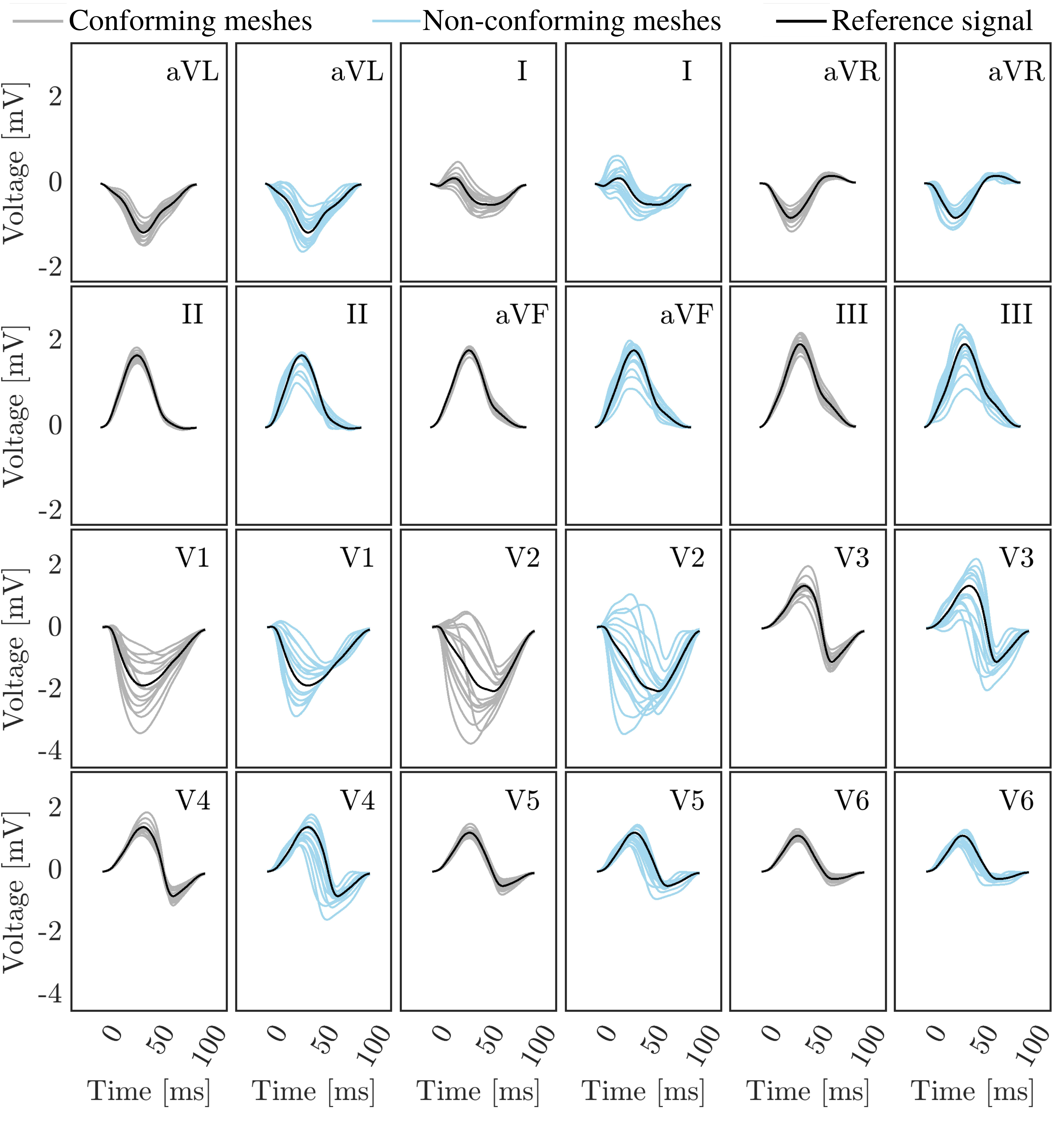} 
	\caption{Qualitatively comparison between the QRS complex obtained displacing the cardiac geometry with big transformations. In grey we show the ECG from conforming heart-torso geometries, in ligh-blue those computed with heart-torso non-conforming geometries. We also draw the ECG signal obtained with the reference heart-torso geometry (black).}
	\label{Fig.variability_to_geometry_movement}
\end{figure}

Lastly, in Fig. \ref{Fig.variability_to_geometry_movement}, we present the ECG traces for all major transformations, demonstrating that our approach is capable of qualitatively reproducing a highly similar pattern across most leads. Furthermore, significant differences between the two groups of signals are primarily characterized by variations in amplitude rather than the shape of the ECG waveform. This highlights the reliability of our approach, particularly in studying the impact of electrophysiological assumptions on a cohort of patients representing a diverse human population, at regarding the produced ECG.

\section{Conclusions}
\label{Sec:conclusions}
In this study, we introduced a segregated and staggered numerical framework designed for the computation of clinical outputs pertaining to cardiac electrophysiology, including body surface potential maps and electrocardiograms. Our model is derived as a simplification of the fully-coupled heart-torso model, which is known for its comprehensive nature but also its computational costs. Our approach employs a one-way coupled scheme between the electrophysiology model and the Laplace model representing the torso as a passive conductor. This allows for less frequent solution of the model in the torso, resulting in lower computational costs. Additionally, we enhance the flexibility of the framework by incorporating interpolation methods at the heart-torso interface, enabling different spatial discretizations within the heart and torso domains.

To assess the reliability of our model, we carried out a comparison with the fully coupled heart-torso model. We examined the body surface potential maps and electrocardiograms both qualitatively and quantitatively, finding only minimal differences in the model outputs. This highlights the robustness of our approach and its ability to accurately capture the relevant features of the fully coupled model.

Next, we considered a thorough testing of our framework by introducing domain non-conformity in both the heart and torso domains. Firstly, we demonstrate that our model enable the use of significantly different spatial discretizations. Exploiting the regularity of the model solution in the torso domain, we emphasize the advantages of reducing the mesh size specifically within the cardiac domain. This reduction does not hamper the clinical information obtained from the model outputs, while significantly decreasing the computational costs associated with the torso model. Our findings indicate that our framework can be considered highly efficient in comparison to classical numerical schemes used to solve the fully-coupled heart-torso model, and provides valuable insights into the general dependency of the torso solution on the model discretization. 

Finally, we demonstrate that our methods can effectively reduce computational time associated with the pre-processing of torso geometries. This is achieved by safely employing our approach when small geometric non-conformities, such as the one resulting from cardiac breathing, exist between the heart and the torso domains. We assess this capability by studying both  BSPMs and ECGs using a manually defined cohort of heart-torso geometries created through the displacement of the cardiac domains with six rigid transformations. Furthermore, we investigate the variability of clinical outputs resulting from significant geometric non-conformities, such as those representing a diverse human population. We examine this variability qualitatively and demonstrate the ability of our approach to reproduce good results, particularly in terms of ECG traces, depending on the specific geometric transformations being considered.

Major limitations in this work are related with the modeling of certain important features of cardiac electrophysiology, such as cellular heterogeneity, which play a role in accurately reproducing the T wave of the ECG. Additionally, the inclusion of the atria in the cardiac geometry, which influences the P wave through atrial depolarization, is also not accounted for in this study. However, it is important to note that these details can be easily incorporated into our framework, as they primarily impact cardiac electrophysiology and do not significantly affect the heart-torso interactions captured in this work.

Furthermore, in this paper we assumes homogeneous tissue within the torso, which allows for easier coarsening of the torso domain by omitting small tissue features. However, tissue heterogeneity of different to organs can be easily included in the torso domain, although it has been demonstrated to have minimal influence on the computation of the outputs and is often disregarded \cite{Keller2010}.

Lastly, it is worth noting that this study does not directly compare the computational efficiency of our framework against other commonly used mathematical models, such as the lead field or the boundary element methods, which are often employed in the computation of electrophysiology clinical outputs. While model comparison will be a potential avenue of exploration in future work, it is important to emphasize that our objective was to demonstrate the potential of our method to accurately reproduce a wide range of diverse cardiac electrical outputs, many of which are not readily reproducible using alternative models.

In conclusion, we have demonstrated that our numerical approach serves as a reliable, efficient, and flexible alternative to the highly detailed heart-torso coupled model. Our specific formulation exhibits flexibility in terms of modeling assumptions, and it paves the way for integrating the electrophysiology-torso model with other aspects of cardiac functioning, such as mechanics and fluid dynamics. These integrations will be addressed in forthcoming publications.

\section*{Formatting of funding sources}
This research has been funded by the European Research Council (ERC) under the European
Union’s Horizon 2020 research and innovation programme, grant agreement No. 740132, iHEART “An
Integrated Heart Model for the simulation of the cardiac function”, P.I. Prof. A. Quarteroni.

%\section*{Sample CRediT author statement}
%\textbf{Elena Zappon:} Conceptualization, Methodology, Software, Formal analysis, Writing - Original Draft. \textbf{Andrea Manzoni:} Writing- Reviewing and Editing. \textbf{Alfio Quarteroni:} Writing- Reviewing and Editing.

\bibliography{Ref}

\appendix
\section{Numerical schemes}
\label{Appendix}
In this appendix we show the discretized numerical formulation of the models analyzed in this work when the mesh are conforming. The same formulation for the one-way coupled problems can then be applied to the non-conforming case by adding the linear interpolation operator $\Pi$ to interpolate $u_e$ from the cardiac interface to the torso one, as described in Subsection \ref{Subsec:framework}.
 
We employ continuous finite elements of either order 2 ($\mathbb{P}_{2}$) - the same can be done using $\mathbb{P}_{1}$ -  and two triangulations $\mathcal{T}_H$ and $\mathcal{T}_T$ over the heart $\Omega_H$ and the torso $\Omega_T$ domains, respectively, made of 3D tetrahedra. On the cardiac domain, we then define 
set $\mathbb{P}_p(K)$ of tensor-products of polynomials with degree smaller or equal to $p$ over a mesh element $K$, and the finite dimensional spaces 
$$X_H^p = \{v \in C^0(\overline{\Omega}_H): v_{|K} \in \mathbb{P}_p(K) ~ \forall K \in \mathcal{T}_H\},$$
and
$V_H^p = X_H^p \cap V_H$,  $V_I^p = X_H^p \cap V_I$, where $V_{H}$, $V_{e,H}$ and $V_I$ are the Sobolev spaces for the corresponding weak formulation \cite{QuarteroniValli2008}, and  $V_{e,H}^p = X_H^p \cap V_{e,H}$, repsectively. % and I denote their dimensions by $N_H = \text{dim}(V_H^p)$, $N_I = \text{dim}(V_I^p)$, and $N_{e,H} = \text{dim}(V_{e,H}^p)$.
 Similarly, for the torso problem we can define the finite dimensional spaces
$V_{T}^p = V_{T} \cap X_T^p$, and  $V_{HT}^p = V_{HT} \cap X_{HT}^p,$  with   
$$ X_T^p = \left\{v \in C^0(\overline{\Omega}_T): v_{|K} \in \mathbb{P}_p(K), p\geq 1, \forall K \in \mathcal{T}_T\right\}$$
and
$$ X_{HT}^p = \left\{v \in C^0(\overline{\Omega_H \cup \Omega_T}): v_{|K} \in \mathbb{P}_p(K), p\geq 1, \forall K \in \mathcal{T}_H \cup \mathcal{T}_H\right\}.$$

Given the time discretization presented in Subsection \ref{Subsec:framework}, here we will use the same time step for both cardiac EP and torso model - i.e. $m=1$ - and we will denote by $\mathbf{a}^i \simeq \mathbf{a}(t_i)$ the generic fully discretized FEM approximation of a scalar, vectorial or tensorial variable $\mathbf{a}(t)$ over $\Omega_H$ or $\Omega_T$. 

\subsection*{Fully-coupled heart-torso model}
The fully-coupled heart-torso model is defined by assuming perfect coupling conditions between the bidomain model \ref{Eq:bidomain_model} and torso model \eqref{Eq:torso}, i.e. by imposing continuity of the solution and of the flux at the domains interface with Dirichlet and Neumann interface conditions, respectively, that with
\begin{equation}
\label{Eq:interface_FCHT}
\begin{cases}
u_e = u_T &\text{on }\Gamma\\
(\mathbf{D}_i \nabla V_m + (\mathbf{D}_i + \mathbf{D}_e)\nabla u_e) \cdot \mathbf{n}_H = -\mathbf{D}_T \nabla u_T \cdot \mathbf{n}_T &\text{on }\Gamma.
\end{cases}
\end{equation}
In this work, we solve such reaction-diffusion in monolithic form, whereas the solution of ionic model \eqref{Eq:ionic_model} is treaten with an IMEX scheme (see \ref{Subsec:framework}) and therefore computed before solving the reaction-diffusion model.
 
Specifically, given the initial conditions
$$ V_m^0 = V_{m,0}, \quad \mathbf{w}^0 = w_0 \quad\text{and }  \mathbf{c}^0 = c_0,$$
at each time step $t_{i+1}, ~n=0,\dots,N_H$ we search for $(V_m^{i+1}, u_e^{i+1}, u_T^{i+1})$ through the following steps:
\begin{enumerate}
	\item Approximate $V_m^{i+1}$ through a first order extrapolation, \emph{i.e.} $V^{i+1}_{m,\text{EXT}} =  V_m^{i-1}.$
	\item Solve the ionic model and compute the ionic current $I^{\text{ion}}(V^{i+1}_{m,\text{EXT}},\mathbf{w}^{i+1},\mathbf{c}^{i+1})$, where $\mathbf{w}^{i+1},\mathbf{c}^{i+1}$ are the solutions of the ionic model at time $t_i$.
	\item Assemble and solve the following system
	\begin{equation}
	\label{Eq:monolithic_discrete_formulation} 
	\begin{cases}
	\begin{split}
	\chi\int_{\Omega_H} \frac{ 3 C_m}{2\Delta t} V_m^{i+1}\varphi ~ &+~  \int_{\Omega_H} \mathbf{D}_i \nabla(V_m^{i+1} + u_e^{i+1})\cdot \nabla \varphi \\ & =\chi\int_{\Omega_H} \frac{C_m}{\Delta t} \left(2 V_m^{i} - \frac{1}{2}V_m^i\right) \varphi \\&\quad- \chi\int_{\Omega_H} I^{\text{ion}}(V^{i+1}_{m,\text{EXT}},\mathbf{w}^{i+1},\mathbf{c}^{i+1}) \varphi\\&\quad+ \int_{\Omega_H}I^{\text{app}}(t_i)\varphi \qquad \qquad \qquad \qquad \:\:\,\forall \varphi \in V_H^p
	\end{split} \\
	\begin{split}\int_{\Omega_H} \mathbf{D}_i \nabla V_m^{i+1} \cdot \nabla \phi ~&+~ \int_{\Omega_h}(\mathbf{D}_i + \mathbf{D}_e)\nabla u_e^{i+1} \cdot \nabla \phi \\& \:\:\, + \int_{\Omega_T}\mathbf{D}_T \nabla u^{i+1}_T \cdot \nabla \phi = 0 \qquad \qquad \forall \phi \in V_{HT}^p,
	\end{split} 
	\end{cases}
	\end{equation}
\end{enumerate}
where $\{\varphi_i\}_{i=1}^{N_H}$ and $\{\phi_i\}_{i=1}^{N_S}$ are the set of basis functions of $V_H^p$ and $V_{HT}^p$, being $N_H$ and $N_S$ the dimensions of $V_H^p$ and $V_{HT}^p$, respectively.

In a practical implementation, system \eqref{Eq:monolithic_discrete_formulation} is assembled  introducing a block matrix with two main blocks: one for the bidomain model and one for the torso model. Homogeneous constraints are used to include Dirichlet interface conditions of \eqref{Eq:interface_FCHT}, enlarging the final system, while Neumann interface conditions of \eqref{Eq:interface_FCHT} are naturally dropped in the equations set. Moreover, an algebraic multigrid (AMG) preconditioner \cite{Stuben2001,Plank2006} is implemented to accelerate and stabilize the solution of system \eqref{Eq:monolithic_discrete_formulation} (see \cite{Vigmond2008,Scacchi2008,GerardoGiorda2009} for other computational techniques related to the solution of the coupled Bidomain-torso system). 

\subsection*{One-way bidomain-torso and pseudo-bidomain-torso model}
The one-way coupling is obtained by assuming an isolated heart from the torso, thus imposing homogeneous Neumann interface conditions \eqref{Eq:Neumann_one_way} on the cardiac EP model - either the bidomain or the pseudo-bidomain one - while prescribing non-homogeneous Dirichlet interface conditions \eqref{Eq:dirichlet_torso} at the interface of the torso domain. As a result, the intra-cardiac quantities $V_m$ and $u_e$ can be computed independently on the torso model, employing the bidomain model \eqref{Eq:bidomain_model} or the pseudo-bidomain model \eqref{Eq:monodomain_model}-\eqref{Eq:pseudo_bid} in a insulated formulation.

The numerical schemes of the one-way coupled models are then obtained by employing the same numerical framework of the FCHT model, that is $\mathbb{P}_2$ FE for the space discretization, and an IMEX treatment of the ionic term of the EP model. Therefore, step 1 and 2 of the numerical schemes in the previous subsection remain unchanged. Step 3 is instead replaced by the solution of the bidomain or pseudo-bidomain model. 

In the former case, the electrophysiological system is
\begin{enumerate}
	\item [3.] for each $t_i$, $n=0,\dots,N_H$, find $(V_m^{i+1},u_e^{i+1}) \in V_H^p \times V_{e,H}^p$ such that
	\begin{equation}
	\begin{cases}
	\begin{split}
	\chi\int_{\Omega_H} \frac{ 3 C_m}{2\Delta t} V_m^{i+1}\varphi ~ &+~  \int_{\Omega_H} \mathbf{D}_i \nabla(V_m^{i+1} + u_e^{i+1})\cdot \nabla \varphi \\ & =\chi\int_{\Omega_H} \frac{C_m}{\Delta t} \left(2 V_m^{i} - \frac{1}{2}V_m^i\right) \varphi \\&\quad- \chi\int_{\Omega_H} I^{\text{ion}}(V^{i+1}_{m,\text{EXT}},\mathbf{w}^{i+1},\mathbf{c}^{i+1}) \varphi\\&\quad+ \int_{\Omega_H}I^{\text{app}}(t_i)\varphi \qquad \qquad \qquad \qquad \:\:\:\, \forall \varphi \in V_H^p
	\end{split} \\
	\begin{split}\int_{\Omega_H} \mathbf{D}_i \nabla \mathbf{V}_m^{i+1} \cdot \nabla \phi ~&+~ \int_{\Omega_h}(\mathbf{D}_i + \mathbf{D}_e)\nabla u_e^{i+1} \cdot \nabla \phi = 0 \qquad \forall \phi \in V_{e,H}^p.
	\end{split} 
	\end{cases}
	\end{equation}
\end{enumerate}

In the latter case, first we solve the monodomain model and then we recover the extra-cellular potential;
\begin{enumerate}
	\item [3.] for each time instant $t_i$, $n=0,\dots,N_H$:
	\begin{enumerate}
		\item [a.] find $V_m^{i+1} \in V_H^p$ such that
		\begin{equation}
		\begin{split}
		\chi\int_{\Omega_H} \frac{ 3 C_m}{2\Delta t} V_m^{i+1}\varphi ~ &+~  \int_{\Omega_H} \mathbf{D}_m \nabla V_m^{i+1} \cdot \nabla \varphi \\ & =\chi\int_{\Omega_H} \frac{C_m}{\Delta t} \left(2 V_m^{i} - \frac{1}{2}V_m^i\right) \varphi \\&\quad- \chi\int_{\Omega_H} I^{\text{ion}}(V^{i+1}_{m,\text{EXT}},\mathbf{w}^{i+1},\mathbf{c}^{i+1}) \varphi\\&\quad+ \int_{\Omega_H}I^{\text{app}}(t_i)\varphi \qquad \qquad \qquad \qquad  \:\:\:\,\forall \varphi \in V_H^p. 
		\end{split} 
		\end{equation}
		\item [b.] given $V_m^{i+1}$, find $u_e^{i+1} \in V_{e,H}^p$ such that
		\begin{equation}
		\int_{\Omega_H}(\mathbf{D}_i + \mathbf{D}_e) \nabla u_e^{i+1} \cdot \nabla \psi = \int_{\Omega_H}\mathbf{D}_i  \nabla V_m^{i+1} \cdot \nabla \psi \qquad \forall\psi\in V_{e,H}^p
		\end{equation}
	\end{enumerate}
\end{enumerate}

A fourth step is added to solve the torso model \eqref{Eq:torso}. 
\begin{enumerate}
	\item [4.] Knowing $u_e^{i+1}$, find $u_T^{i+1} \in V_{HT}^p$ such that
	\begin{equation}
	\label{Eq:torso_uncoupled}
	\int_{\Omega_T} \mathbf{D}_T \nabla u_T^{i+1} \cdot \nabla \phi = 0 \qquad\forall  \phi\in V_{HT}^p
	\end{equation}
	with $u_T^{i+1} = u_e^{i+1} \text{ on }\Gamma$.
\end{enumerate}
See Fig. \ref{Fig:scheme_num} for a schematic representation of the uncoupled model solution with either Bidomain or Pseudo-Bidomain model. As for the FCHT model, we employ AMG preconditions to accelerate and stabilize the solution of the computed algebraic systems.

\end{document}